\documentclass{amsart}
\usepackage[margin=3cm]{geometry}
\usepackage{graphicx}
\usepackage{bm}
\usepackage[hidelinks]{hyperref} 
\usepackage[english]{babel}
\usepackage{commath}
\usepackage{latexsym}
\usepackage{amssymb}
\usepackage{amsmath}
\usepackage[foot]{amsaddr}
\usepackage{amsfonts}
\usepackage{graphics}
\usepackage{mathrsfs}
\usepackage{mathtools}
\usepackage{amsopn}
\usepackage{xcolor}
\usepackage{verbatim}

\makeatletter
\edef\orig@output{\the\output}
\output{\setbox\@cclv\vbox{\unvbox\@cclv\vspace{0pt plus 20pt}}\orig@output}
\makeatother

\DeclareMathOperator{\dive}{div}

\DeclareMathOperator{\rb}{\text{rb}}

\newcommand{\lbar}[1]{\overline{#1}}

\newcommand{\gr}[1]{\bm{#1}}
\newcommand{\grs}[1]{\boldsymbol{#1}}

\usepackage{multicol} 
\usepackage{xpatch}
\usepackage{nomencl} 
\renewcommand*{\nompreamble}{\footnotesize}
\makenomenclature
\RequirePackage{ifthen}
\let\oldref\ref
\renewcommand{\ref}[1]{(\oldref{#1})}

\renewcommand*\nompreamble{\begin{multicols}{2}}
\renewcommand*\nompostamble{\end{multicols}}

\xpatchcmd{\thenomenclature}{%
  \section*{\nomname}
}{
}{\typeout{Success}}{\typeout{Failure}}

\graphicspath{{./img/}}

\usepackage{ifthen}
\renewcommand{\nomgroup}[1]{%
  \ifthenelse{\equal{#1}{A}}{\item[\textbf{Abbreviations}]}{%
    \ifthenelse{\equal{#1}{G}}{\item[\textbf{Symbols}]}{%
      \ifthenelse{\equal{#1}{C}}{\item[\textbf{Abbreviations}]}{%
        \ifthenelse{\equal{#1}{S}}{\item[\textbf{Subscripts}]}{%
          \ifthenelse{\equal{#1}{Z}}{\item[\textbf{Mathematical Symbols}]}{}
        }
      }
    }
  }
}

\begin{document}

\title{A reduced order variational multiscale approach for turbulent flows}

\author{Giovanni Stabile$^1$}
\author{Francesco Ballarin$^1$}
\author{Giacomo Zuccarino$^1$}
\author{Gianluigi Rozza$^1$}
\address{$^1$ mathLab, Mathematics Area, SISSA, via Bonomea 265, I-34136 Trieste, Italy}

\maketitle

\begin{abstract}
The purpose of this work is to present different reduced order model strategies {starting from full order simulations stabilized using a residual based Variational MultiScale (VMS) approach. The focus is on flows with moderately high Reynolds numbers. The Reduced Order Models (ROMs) presented in this manuscript are based on a POD-Galerkin approach.}
Two different reduced order models are presented, which differ on the stabilization used during the Galerkin projection.
In the first case the VMS stabilization method is used at both the full order and the reduced order level. In the second case, the VMS stabilization is used only at the full order level, while the projection of the standard Navier-Stokes equations is performed instead at the reduced order level. 
The former method is denoted as \emph{consistent ROM}, while the latter is named \emph{non-consistent ROM}, in order to underline the different choices made at the two levels.
Particular attention is also devoted to the role of inf-sup stabilization by means of supremizers in ROMs based on a VMS formulation. Finally, the developed methods are tested on a numerical benchmark.
\end{abstract}

\section{Introduction}\label{sec:intro}

During the last decades, the increase of computational resources and the great improvement into numerical methods to solve partial differential equations has allowed to broaden the applicability of computational methods, originally devised in an academic environment, to industrial problems. Indeed, nowadays numerical simulation is a key ingredient in several engineering and applied sciences fields such as aeronautical, naval, civil and mechanical engineering, life sciences and medicine where computational methods are used to solve mechanics and fluid dynamics problems.

However, especially for fluid dynamics, there are still certain situations where the direct numerical simulations of the governing equations using standard discretization techniques (Finite Elements, Finite Volumes, Finite Differences, Spectral Element Methods) may become unaffordable. Such situations occurs when a large number of different system configurations are in need of being tested (shape optimization, uncertainty quantification, inverse problems) or an extremely reduced computational cost is required (real-time control). 

A viable approach to solve this kind of problems, at a far lower computational cost, is given by reduced order modelling techniques \cite{bennerParSys,ChinestaEnc2017,libroRozza,quarteroniRB2016}. Reduced order models are based on the assumption that the system response is sufficiently smooth with respect to the input parameters. Given this assumption, the solution space can be restricted to a lower dimensional space given by properly chosen basis functions. This work is based on a POD-Galerkin approach which consists into the extraction of the most energetic modes representing the system dynamics and into a Galerkin projection onto the space spanned by these modes. POD-Galerkin methods have been proposed in several works starting from different full order discretization techniques and focusing the attention on different aspects. We mention \cite{deane1991,Iollo1998,ITO1998403,noack1994,peterson1989} among the first attempts to create low-dimensional models for fluid flows starting from POD basis functions.

However, it is well known that standard POD-Galerkin models are prone to possible stability issues. The type of instabilities which are analyzed an treated in literature can be classified into different classes depending on the source of instability: (i) in \cite{Akhtar2009,Bergmann2009516,taddei2017,Iollo2000,Sirisup2005218,giovannisaddam2017} the attention is focused on instability observed in transient problems during long time integration; (ii) \cite{Shafqat,Ballarin2015,Caiazzo2014598,Gerner2012,rozza2007stability,Stabile2018} are devoted to accurate pressure recovery and inf-sup instability, while (iii) in \cite{ChaconDelgadoGomezBallarinRozza2017,iliescu2013variational,Iliescu_vms,wang_turb} the instability due to advection dominated phenomena {(see also \cite{Giere2015} for a similar approach on convection-diffusion problems)} and transition to a turbulent regime are discussed. 
{Focusing more the attention on instabilitities due to turbulent and advection dominated flows a promising approach is to employ subgrid-scale closure models. In \cite{iliescu2013variational,Iliescu_vms} the authors developed a variational multiscale POD reduced order model and presented also a detailed analysis on error estimates. In the first work a VMS closure model is used to stabilize convection-dominated convection-diffusion-reaction equations and in the second case it has been extended to Navier-Stokes equations in a turbulent setting. In \cite{wang_turb} a dynamic subgrid-scale closure model is presented and compared with a variational multiscale method. 
{This work follows what already done by some of the authors of this article in \cite{PACCIARINI20141} where a SUPG stabilization is used for advection-diffusion PDEs and is also inspired by \cite{Giere2015} that follows a similar line of thought as the current approach but using a SUPG stabilization for convection-dominated convection-diffusion-reaction equations.}

The present contribution, starting from what done in the aforementioned references, aims to deal with instabilities of type (ii) and (iii) and to study the relationship between the stabilization techniques adopted for inf-sup pressure instabilities and the ones adopted for convection-dominated/turbulent problems. Differently to what has been done in \cite{iliescu2013variational,Iliescu_vms,wang_turb} we employ {the same} variational multiscale stabilization technique at both the full order and the reduced order level. The reason for this choice is given by the fact that we want to operate in conditions where a DNS simulation becomes not realizable and therefore, a stabilization is also needed during the full order simulations. Another difference respect to previous published sources is that the reduced reduced order model presented in this work accounts also for the pressure contribution. The methodology could be then employed to study also more general cases where the pressure term cannot be neglected \cite{noack1994}. The main novelty of the present work consists into the development of a reduced order model starting from VMS stabilized full order model and into a discussion concerning the inf-sup stability of the resulting ROM.}

It is our belief that literature on closure models for the approximation of turbulent flows using standard discretization techniques presents a large variety of options (see e.g. \cite{davidson2004,pope2001,sagaut2006}) which are waiting to be incorporated into a reduced order framework.
Among the most studied strategies we mention the Reynolds-Averaged Navier-Stokes (RANS) models \cite{pope2001}, the Large-Eddy Simulation (LES) models \cite{sagaut2006} and the Variational Multiscale (VMS) method \cite{Hughes1996,Hughes1998}.
In the development of a reduced order strategy, one must take into account that the choice of the closure model often depends on the specific application at hand. Indeed, closure models based on eddy viscosity (e.g. $k-\omega$, $k-\epsilon$, Spalart-Allmaras), mainly tailored for finite volumes discretization, seem to be particularly promising to treat developed turbulent flows with high Reynolds number in the context of industrial applications. A reduced order strategy for such applications is going to be actively pursued by the authors of this manuscript in \cite{HijaziAliStabileBallarinRozza2018,Lorenzi2016,saddam2018}. Here, however, {we focus on closure models} adopted for more moderate Reynolds numbers, based on a finite element discretization with VMS stabilization and it fits well moderate turbulent patterns.
The proposed methodology can still be applied to several relevant cases, such as those arising in biomedical applications. Furthermore, it builds on a stronger theoretical framework in the reduced order modelling community, which has recently seen first contributions for what concerns certification of a posteriori error estimators in turbulent flows \cite{ChaconDelgadoGomezBallarinRozza2017}. 

The paper is organized as follows: the mathematical formulation of the problem is set in \autoref{sec:math_form}; the VMS finite element discretization used to numerically solve the mathematical problem is reported in \autoref{sec:FOM}. \autoref{sec:ROM} contains the formulation of the Reduced Order Model, and some consideration concerning its stabilization. The proposed methodology is tested on a numerical example in \autoref{sec:numer_example}, and finally some conclusions are drawn and perspectives for future developments are given in \autoref{sec:concl}.

\section{Mathematical Formulation of the Physical Problem}\label{sec:math_form}
The mathematical problem on which this work is focused is given by the incompressible Navier-Stokes equations. Given a space-time domain $Q=\Omega\times[0,T]\subset\mathbb{R}^2\times\mathbb{R}^+$ the initial boundary-value problem consists in solving the following equations for the velocity $\bm{u}:Q\rightarrow\mathbb{R}^2$ and the pressure $p:Q\rightarrow\mathbb{R}$ such that:
\begin{equation}\label{eq:navstokes}
\begin{cases}
\partial_t\bm{u}+ \dive (\bm{u} \otimes \bm{u})- \dive(2 \nu \nabla^s \bm{u}) + \nabla p = \bm{f} &\mbox{ in } Q,\\
\dive(\bm{u})=\bm{0} &\mbox{ in } Q,\\
\bm{u} = \bm{g} &\mbox{ on } \Gamma_{\text{in}} \times [0,T],\\
\bm{u} = \bm{0} &\mbox{ on } \Gamma_{\text{wall}} \times [0,T],\\
(2 \nu\nabla^s \bm{u} - p\bm{I})\bm{n} = \bm{0} &\mbox{ on } \Gamma_{\text{out}} \times [0,T],\\ 
\bm{u}=\bm{u}_0 &\mbox{ in } \Omega\times\{0\},\\            
\end{cases}
\end{equation}
where $\Gamma = \Gamma_{\text{in}} \cup \Gamma_{\text{wall}} \cup \Gamma_{\text{out}}$ is the boundary of $\Omega$, composed by three non-overlapping regions $\Gamma_{\text{in}}$, $\Gamma_{\text{wall}}$ and $\Gamma_{\text{out}}$ that indicate, respectively, inlet boundary, physical walls, and outlet boundary.
Furthermore, $\nu \in \mathbb{R}^+$ is the kinematic viscosity, $\gr{f}:Q\rightarrow\mathbb{R}^2$ is the given body force (per unit volume), $\gr{g}:\Gamma_{\text{in}}\times[0,T]\rightarrow\mathbb{R}^2$ is the inlet velocity profile, and $\bm{u}_0: \Omega \to \mathbb{R}^2$ denotes the initial condition for the velocity.
The symbol $\otimes$ denotes the tensor product of two vectors in $\mathbb{R}^2$, i.e. in components $\left(\gr{u}\otimes\gr{v}\right)_{ij}=u_i v_j$, while $\nabla^s\gr{u}$ denotes the symmetric part of the velocity gradient, i.e. $\nabla^{s}\gr{u}=\frac{\nabla\gr{u}+\nabla\gr{u}^{T}}{2}$.

The weak formulation of \eqref{eq:navstokes} is obtained by multiplying the momentum equation for a test function $\gr{v}\in\gr{V} = \{\gr{v} \in \gr{H}^1(\Omega): \gr{v} = \gr{0}\text{ on }\Gamma_{\text{in}} \cup \Gamma_{\text{wall}}\}$ and the continuity equation for a test function $q\in Q = L^2(\Omega)$, and with integration by parts over the space domain $\Omega$.
The weak problem for \eqref{eq:navstokes} reads as follows:

\noindent\emph{Find $\gr{u}\in L^2([0,T],\gr{H}^1(\Omega))$ and $p\in L^2([0,T],Q)$ such that $\gr{u}=\gr{g}$ on $\Gamma_{\text{in}}\times[0,T]$ and $\gr{u}=\gr{0}$ on $\Gamma_{\text{wall}}\times[0,T]$, and the equations
\begin{equation}
\label{NS weak}
\begin{cases}
(\partial_t\gr{u},\gr{v})_{\Omega}-(\gr{u}\otimes\gr{u},\nabla \gr{v})_{\Omega}+2\nu(\nabla^s\gr{u},\nabla^s\gr{v})_{\Omega}-(p,\dive(\gr{v}))_{\Omega}=(\gr{f},\gr{v})_{\Omega},\\
(q,\dive(\gr{u}))_{\Omega}=0,\\
\gr{u}(0)=\gr{u}_0,
\end{cases}
\end{equation}
hold for all $\gr{v}\in\gr{V}$ and for all $q\in Q$.
}


\section{The Full Order Model}\label{sec:FOM}
In this section the full order discretization used through the manuscript is presented. In \autoref{subsec:FEM_GAL} the standard Galerkin finite element approximation is introduced, while \autoref{subsec:FEM_VMS} presents the Residual Based Variational Multiscale concepts used to achieve stable simulations.  

\subsection{The Standard Galerkin Finite Element Approximation}\label{subsec:FEM_GAL}
Given the continuous formulation presented in \autoref{sec:math_form} here we introduce the discrete formulation and the variational multiscale finite element approximation.
We consider a finite decomposition (mesh) $\mathcal{T}_h$ of $\lbar{\Omega}=\bigcup_{K\in \mathcal{T}_h} K,\;$ composed of triangular cells $K$. 
Two finite element (FE) spaces $\gr{V}_h\subset\gr{V}$ and $Q_h\subset Q$ are employed to discretize \eqref{NS weak} on $\mathcal{T}_h$.
Specifically, $\gr{V}_h$ and $Q_h$ are spaces of continuous locally polynomial functions of degree $2$ and $1$, respectively. The FE formulation on $(\gr{V}_h,Q_h)$, known as Taylor-Hood FE, satisfies the so called \emph{inf-sup} condition \cite{boffi_mixed}
\begin{equation}
\label{inf-sup}
\inf_{q_h\in Q_h\setminus\{0\}}\sup_{\gr{v}_h\in\gr{V}_h\setminus\{\gr{0}\}}\frac{(\dive(\gr{v}_h),q_h)}{\norm{\gr{v}_h}_{\gr{V}_h}\norm{q_h}_{Q_h}}\geq\beta>0.
\end{equation}
We will come back to the topic of inf-sup stability for what concerns the reduced order model. For the sake of notation in the definition of the reduction process, let us define two sets $\{\grs{\varphi}_j\}_{j=1}^{M_h}$, $\{\psi_l\}_{l=1}^{K_h}$ of Lagrangian basis functions of $\gr{V}_h$ and $Q_h$, respectively, being $M_h=\text{dim}(\gr{V}_h)$ and $K_h = \text{dim}(Q_h)$.

The semi-discrete Galerkin-FE approximation of \eqref{NS weak} reads as follows: find $\gr{u}_h\in L^2([0,T],\gr{V}_h)$ and $p_h\in L^2([0,T],Q_h)$ such that $\gr{u}_h=\gr{g}_h$ on $\Gamma_{\text{in}}\times[0,T]$ and $\gr{u}_h=\gr{0}$ on $\Gamma_{\text{wall}}\times[0,T]$, and the equations
\begin{equation}
\label{NS weak discrete}
\left\lbrace
\begin{split}
\left(\partial_t\gr{u}_h,\gr{v}_h\right)_{\Omega}&-\left(\gr{u}_h\otimes\gr{u}_h,\nabla \gr{v}_h\right)_{\Omega}+2\nu\left(\nabla^s\gr{u}_h,\nabla^s\gr{v}_h\right)_{\Omega}+\\
&-\left(p_h,\dive(\gr{v}_h)\right)_{\Omega}+\left(q_h,\dive(\gr{u}_h)\right)_{\Omega}=\left(\gr{f},\gr{v}_h\right)_{\Omega},\\
\gr{u}_h(0)&=\gr{u}_{0,h},
\end{split}
\right.
\end{equation}
hold for all $\gr{v}_h$ in $\gr{V}_h $ and for all $q_h$ in $Q_h$, being $\gr{g}_h$ and $\gr{u}_{0,h}$ suitable interpolations of $\gr{g}$ and $\gr{u}_{0}$ on the mesh $\mathcal{T}_h$.

Finally, time discretization is carried out by taking $N_T+1$ time instants $\{t_k\}_{k=0}^{N_T}$ in $[0,T]$, chosen such that $\Delta t=t_k-t_{k-1}$ is constant. For simplicity a backward Euler time stepping scheme is used in this manuscript. Since such time discretization is very standard (see e.g. \cite{quarteroni2008numerical}), and seamlessly applies at the FE, VMS and reduced order levels, in the following we will always refer to the semi-discrete formulation \eqref{NS weak discrete} omitting any further presentation of the time discretization.

\subsection{The Residual Based Variational Multiscale Formulation}\label{subsec:FEM_VMS}
In case of turbulent flows the standard Galerkin-FE approximations may fail to accurately model the physical phenomena unless a very refined mesh is employed. As such refinement is usually unaffordable for Reynolds numbers in the order of a few thousands, we resort instead to a stabilization based on a variational multiscale method.
In this work we rely on the Residual Based variational multiscale method, as presented in \cite{Bazilevs2007173}. 
The Residual Based VMS was first introduced in \cite{Hughes1995}. Further theoretical development that led to the application in Computational Fluid Dynamics (CFD) were given in \cite{Hughes1998,Hughes1996}. 
A residual based VMS method has been used in CFD applications for a large variety of problems (see e.g. \cite{Bazilevs2007173,codina2017variational,forti2015semi,hughes2004multiscale,masud2011heterogeneous}), and as a result such VMS method has proven to be an efficient and accurate turbulence model for applications characterized by moderately high Reynolds number flows.
The method is a two-scales method, even though other VMS variants with a larger number of scales are available in literature (e.g. the Residual Free Bubble VMS \cite{Brezzi1992}, the Projection Based VMS \cite{rebollo2015numerical} and the Algebraic VMS \cite{Gravemeier2010}). 

Following the presentation in \cite{Bazilevs2007173}, we orthogonally decompose the spaces $\gr{V}$ and $Q$ as
\begin{equation}
\gr{V}=\gr{V}_h\oplus\gr{V}',\;Q=Q_h \oplus Q',
\end{equation}
using the elliptic projector $(\nabla^s\cdot, \nabla^s\cdot)_{\Omega}$ and the $L^2$ projector $(\cdot, \cdot)_{\Omega}$, respectively, as inner products.
The spaces $\gr{V}_h$ and $Q_h$ represent the coarse (resolved) scale, while $\gr{V}'$ and $Q'$ are modelling the fine (unresolved) scale.

Therefore, for any $t\in[0,T]$ we decompose the trial functions in \eqref{NS weak discrete} as
\begin{align*}
\gr{u}(t)&=\gr{u}_h(t)+\gr{u}'(t)&&\text{for }\gr{u}_h(t)\in\gr{V}_h\text{ and }\gr{u}'(t)\in\gr{V}',\\
p(t)&=p_h(t)+p'(t)&&\text{for }p_h(t)\in Q_h\text{ and }p'(t)\in Q',
\end{align*}
assuming the fine scale velocity functions to satisfy homogeneous boundary conditions on the Dirichlet boundaries,
\begin{equation*}
\gr{u}'(t)=0\text{ on }\Gamma_{\text{in}}\cup\Gamma_{\text{wall}},
\end{equation*}
so that the coarse scale velocity functions are such that
\begin{equation*}
\gr{u}_h(t)=\gr{g}_h(t)\;\text{ on }\Gamma_{\text{in}},\quad
\gr{u}_h(t)=\gr{0}\text{ on }\Gamma_{\text{wall}},
\end{equation*}
i.e. they fulfill the physical boundary conditions.

Introducing such expressions in \eqref{NS weak} with coarse scale test functions $\gr{v}_h\in\gr{V}_h$ and $q_h\in Q_h$ we obtain the so-called \emph{coarse scale equation}
\begin{equation}
\label{NS coarse}
\left\lbrace
\begin{split}
(\partial_t\gr{u}_h,\gr{v}_h)_{\Omega}&-(\gr{u}_h\otimes \gr{u}_h,\nabla \gr{v}_h)_{\Omega}+2\nu(\nabla^s \gr{u}_h,\nabla^s \gr{v}_h)_{\Omega}+\\
&-(p_h,\dive(\gr{v}_h))_{\Omega}+(q_h,\dive(\gr{u}_h))_{\Omega}
\\
&-(p',\dive(\gr{v}_h))_{\Omega}-(\nabla q_h,\gr{u}')_{\Omega}+\\
&-(\gr{u}'\otimes\gr{u}_h,\nabla \gr{v}_h)_{\Omega}-(\gr{u}_h\otimes\gr{u}',\nabla \gr{v}_h)_{\Omega}+\\
&-(\gr{u}'\otimes\gr{u}',\nabla \gr{v}_h)_{\Omega}\\
&=(\gr{f},\gr{v}_h)_{\Omega},\\
\gr{u}_h(0)&=\gr{u}_{0,h} .
\end{split}
\right. 
\end{equation}
A few terms omitted in \eqref{NS coarse} drop due to the direct sum orthogonal decomposition, e.g.
\begin{equation*}
\left(\nabla^s\gr{u}',\nabla^s\gr{v}_h\right)_{\Omega}=0,
\end{equation*}
while others, such as in particular
\begin{equation*}
(\partial_t\gr{u}',\gr{v}_h)_{\Omega},
\end{equation*}
are dropped as turbulence modelling assumptions \cite{Bazilevs2007173}.

Following the philosophy of the LES turbulence models the fine scale solutions are not solved explicitly, rather modeled as follows:
\begin{align*}
\gr{u}'&=-\tau_M\gr{r}_M(\gr{u}_h,p_h),\\
p'&=-\tau_Cr_C(\gr{u}_h),
\end{align*}
where such equality is to be understood as restricted to each cell in $\mathcal{T}_h$.

The structure of the stabilization terms introduced in \cite{Bazilevs2007173} and adopted in this work is the following:
\begin{equation}
\label{fine scale explicit}
\begin{cases}
\begin{split}
\gr{r}_M(\gr{u}_h,p_h)&=\partial_t\gr{u}_h+\dive(\gr{u}_h\otimes\gr{u}_h)-\dive(2\nu\nabla^{s}\gr{u}_h)+\nabla p_h-\gr{f},
\\
r_C(\gr{u}_h)&=\dive(\gr{u}_h),
\\
\tau_M&=\tau_M(\gr{u}_h)=\left(\frac{4}{{\Delta t}^2}+\gr{u}_h\cdot G \gr{u}_h+C_{\text{inv}}\nu^2 G : G \right)^{-\frac{1}{2}},
\\
\tau_C&=\tau_C(\gr{u}_h)=(\tau_M\gr{g}\cdot\gr{g})^{-1}.
\\
\end{split}
\end{cases}
\end{equation}
being $\gr{r}_M(\gr{u}_h,p_h)$ and $r_C(\gr{u}_h)$ residuals of the strong form of momentum and continuity equations, respectively, and $\tau_M$ and $\tau_C$ corresponding stabilization coefficients. Furthermore, on each cell $K$,
\begin{equation}
\left\lbrace
\begin{split}
G_{ij}&=\sum_{k=1}^2\frac{\partial\xi_k}{\partial x_i}\frac{\partial\xi_k}{\partial x_j},\quad
 G : G =\sum_{i,j=1}^2 G _{ij} G _{ij},
\\
g_i&=\sum_{j=1}^2\frac{\partial\xi_j}{\partial x_i},\quad
\gr{g}\cdot\gr{g}=\sum_{i=1}^2\gr{g}_i\gr{g}_i,
\end{split}
\right.
\end{equation}
where $\gr{\xi}=(\xi_1,\xi_2)$ is the inverse of the affine mapping between the reference simplex $\hat{K}$ and the triangulation element $K$, and $C_{\text{inv}}$ is the inverse inequality constant \cite{Bazilevs2007173,quarteroni2008numerical}. 

In order to proceed with an algebraic formulation of the resulting semi-discretization we expand $\gr{u}_h(t)$ and $p_h(t)$ in coordinates with respect to the Lagrangian basis on $\gr{W}_h$ and on $Q_h$. Collecting the resulting coordinates in $\mathbf{u}_h\left(t\right)=\left(u_1\left(t\right),\dots,u_{M_h}\left(t\right)\right)^T$ and $\mathbf{p}_h\left(t\right)=\left(p_1\left(t\right),\dots,p_{K_h}\left(t\right)\right)^T$, we obtain the nonlinear dynamical system
\begin{equation}
\begin{cases}
\begin{split}
M\frac{d\mathbf{u}_h\left(t\right)}{dt}&+A\mathbf{u}_h\left(t\right)+C\left(\mathbf{u}_h\left(t\right)\right)\mathbf{u}_h\left(t\right)+D\left(\mathbf{u}_h\left(t\right),\mathbf{p}_h\left(t\right)\right)
+B^T\mathbf{p}_h\left(t\right)
=\gr{F}\left(t\right),\\
B\mathbf{u}_h\left(t\right)&=E\left(\mathbf{u}_h\left(t\right),\mathbf{p}_h\left(t\right)\right),
\end{split}
\end{cases}
\end{equation}
being\footnote{With a slight abuse, in the notation we indicate  for the nonlinear terms their dependence on the coefficient vectors $(\mathbf{u}_h\left(t\right), \mathbf{p}_h\left(t\right))$ on the left-hand side, while the actual expression on the right-hand side will be provided in terms of the corresponding functions $(\gr{u}_h(t), p_h(t))$.}
\begin{align*}
M_{ij}&=\left(\grs{\varphi}_j,\grs{\varphi}_i\right)_{\Omega}&&\text{ the mass matrix},\\
A_{ij}&=\left(\nabla^s\grs{\varphi}_j,\nabla^s\grs{\varphi}_i\right)_{\Omega}&&\text{ the discretized diffusion operator},\\
C\left(\mathbf{u}_h\right)_{ij}&=-\left(\grs{\varphi}_j\otimes\gr{u}_h,\nabla\grs{\varphi}_i\right)_{\Omega}&&\text{ the discretized convection operator},\\
B_{li}&=-\left(\psi_l,\dive\left(\grs{\varphi}_i\right)\right)_{\Omega}&&\text{ the discretized version of the pressure/divergence operator, and}\\
F_i\left(t\right)&=\left(\gr{f},\grs{\varphi}_i\right)_{\Omega}&&\text{ the discretized forcing term},
\end{align*}
for $i,j = 1, \hdots, M_h$ and $l = 1, \hdots, K_h$.
The remaining terms are due to the stabilization, and read: 
\begin{equation*}
\begin{split}
D\left(\mathbf{u}_h,\mathbf{p}_h\right)_i&=
\sum_{K}\left(\left(\nabla\grs{\varphi}_i\right)\gr{u}_h,\tau_M\gr{r}_M\left(\gr{u}_h,p_h\right)\right)_{K}
+\sum_{K}\left(\left(\nabla\grs{\varphi}_i\right)^T\gr{u}_h,\tau_M\gr{r}_M\left(\gr{u}_h,p_h\right)\right)_{K}
\\&
-\sum_{K}\left(\grs{\varphi}_i,\tau_M\gr{r}_M\left(\gr{u}_h,p_h\right)\otimes\tau_M\gr{r}_M\left(\gr{u}_h,p_h\right)\right)_{K}
+\sum_{K}\left(\tau_Cr_C\left(\gr{u}_h\right),\dive\left(\grs{\varphi}_i\right)\right)_{K},
\\
E\left(\mathbf{u}_h,\mathbf{p}_h\right)_l&=\sum_{K}\left(\nabla{\psi}_l,\tau_M\gr{r}_M\left(\gr{u}_h,p_h\right)\right)_{K}.
\end{split}
\end{equation*}

\section{The Reduced Order Models}\label{sec:ROM}
In this section we introduce the two reduced order models (ROMs) employed in this manuscript. We seek a \emph{construction-evaluation} paradigm\footnote{The reason for this terminology, rather than the more widespread \emph{offline-online} decoupling, will be clear in \autoref{subsec:Gproj}.}, as customary in the reduced order modelling community \cite{libroRozza}. A reduced basis is built by means of proper orthogonal decomposition during the construction stage, as detailed in \autoref{subsec:POD}. Afterwards, the ROM is evaluated by a Galerkin projection onto the reduced basis (\autoref{subsec:Gproj}). The two proposed ROMs will only differ during the evaluation phase.

Before laying out the details, let us underline how the proposed ROMs differ in terms of construction and evaluation when compared to previous works \cite{Bergmann2009516,wang_turb}.
The main difference is that in \cite{Bergmann2009516,wang_turb} full order snapshots employed during the construction of the reduced basis are generated using a direct numerical simulation (DNS) approach, while VMS stabilization is used as a closure model for the ROM to account for the contribution of the discarded modes. In the present case we aim to obtain snapshots for values of Reynolds numbers for which a DNS strategy becomes unfeasible. Thus, we need to resort to a VMS stabilization during the construction phase itself. Our interest is then to explore the relationship (in particular, in terms of stability) between a stabilized full order model and the consequent reduced order model built on top of it.

{Finally, we highlight that the accuracy of the ROM will always be measured with respect to the underlying FOM model (i.e. a Navier-Stokes model with VMS stabilization), and not a DNS result. The target of the ROM will thus be to recover as faithfully as possible the FOM solution, regardless of the fact that the FOM solution itself is under resolved (when compared to a DNS).}

\subsection{Construction of the reduced basis through proper orthogonal decomposition}\label{subsec:POD}
The construction of the reduced basis is carried out by means of a Proper Orthogonal Decomposition (POD) \cite{deane1991,ITO1998403,libroRozza,peterson1989}, as follows. Velocity and pressure snapshots, denoted by $\{\gr{u}_h(t_i)\}_{i=1}^{N_T}$ and $\{p_h(t_i)\}_{i=1}^{N_T}$, are obtained solving \eqref{NS coarse}. Then, the following correlation matrices $\Sigma^{\gr{u}}$ and $\Sigma^{p}$ are assembled, as follows:
\begin{equation}
\Sigma^{\gr{u}}_{ij} = \left(\gr{u}_h(t_i) - \langle\gr{u}_h\rangle,\gr{u}_h(t_j) - \langle\gr{u}_h\rangle\right)_{1,\Omega}, \quad \Sigma^{p}_{ij} = \left({p}_h(t_i),{p}_h(t_j)\right)_{\Omega}, \quad i, j = 1, \hdots, N_T,
\end{equation}
where $(\cdot, \cdot)_\omega$ denotes the $L^2(\omega)$ inner product, while $(\cdot, \cdot)_{1, \omega} = (\nabla\cdot, \nabla\cdot)_\omega$ denotes the $\gr{H}^1(\omega)$ inner product. The corresponding induced norms will be denoted similarly as $\lVert\cdot\rVert_\omega$ and $\lVert\cdot\rVert_{1, \omega}$, being either $\omega = \Omega$ for spatial integration or (in the $L^2$ case) $\omega = [0, T]$ for time integration. Furthermore, assuming the inlet profile $\gr{g}$ to be time-independent\footnote{Indeed, non-stationary solutions can be obtained even for time-invariant input data in case of flows with large Reynolds numbers.}, the time-averaged velocity
\begin{equation}
\langle\gr{u}_h\rangle=\frac{1}{N_T}\sum_{i=1}^{N_T}\gr{u}_h(t_i) ,
\end{equation}
is subtracted from the velocity snapshots in order to for $\gr{u}_h(t_i) - \langle\gr{u}_h\rangle$, $i = 1, \hdots, N_T$, to be zero on each Dirichlet boundary. After carrying out an eigendecomposition, we denote by
\begin{equation*}
\lambda_{1}^{\gr{u}}\geq\lambda_{2}^{\gr{u}}\geq\dots\geq\lambda_{N_T}^{\gr{u}}\geq0,\qquad
\lambda_{1}^p\geq\lambda_{2}^p\geq\dots\geq\lambda_{N_T}^{p}\geq0,
\end{equation*}
the resulting eigenvalues, and by 
\begin{equation*}
\lbrace\widehat{\gr{v}}_n\rbrace_{n=1}^{N_T},\qquad \lbrace\widehat{q}_n\rbrace_{n=1}^{N_T},
\end{equation*}
the corresponding eigenfunctions. The reduced basis functions are then obtained as
\begin{equation*}
\gr{\varphi}_n^{\text{rb}}=\sum_{i=1}^{N_T}\left(\gr{\hat{v}}_n\right)_i\gr{u}_h(t_i),
\qquad
\psi_n^{\text{rb}}=\sum_{i=1}^{N_T}\left(\hat{q}_n\right)_i p_h(t_i),
\qquad
n = 1, \hdots, N_T,
\end{equation*}
and possibly normalized in their respective spaces.

It is well known in the reduced basis community that spaces built from (a truncation of) $\{\gr{\varphi}_n^{\text{rb}}\}$ and $\{\psi_n^{\text{rb}}\}$ might not satisfy an inf-sup condition. It is thus customary to enrich the reduced velocity space with the so-called pressure \emph{supremizers} \cite{libroRozza,rozza2007stability,quarteroniRB2016}. Such an approach entails (for each $i = 1, \hdots, N_T$) the computation of a function $\gr{s}_h(t_i)\in\gr{V}_h$ such that
\begin{equation*}
\left(\nabla \gr{s}_h(t_i),\nabla\gr{v}_h\right)_{\Omega}=\left(\nabla p_h(t_i),\gr{v}_h\right)_{\Omega}
\end{equation*}
holds for all $\gr{v}_h\in\gr{V}_h$.
The strategy that we adopt is the \emph{approximate supremizer enrichment} as in e.g. \cite{Ballarin2015}. This approach consists in applying a POD compression to the resulting supremizer snapshots. Following the same procedure as above, we obtain a set $\lbrace\gr{\phi}_n^{\rb}\rbrace_{n=1}^{N_T}$ of supremizer basis functions.

Afterwards, the final result of the construction stage is the definition of the reduced basis spaces\footnote{Numerical results in \autoref{sec:numer_example} will compare the cases with and without supremizer enrichment. It is understood that in the latter case the reduced velocity space is defined as
\begin{equation*}
\gr{V}_{\rb} = \text{span}\{\gr{g}^{\rb}, \gr{\varphi}_1^{\text{rb}}, \hdots, \gr{\varphi}_N^{\text{rb}}\}.
\end{equation*}
We will not cover the latter case in \autoref{subsec:Gproj}, as the required changes to the formulation of the Galerkin projection are very trivial.
}
\begin{equation*}
\gr{V}_{\rb} = \text{span}\{\gr{g}^{\rb}, \gr{\varphi}_1^{\text{rb}}, \hdots, \gr{\varphi}_N^{\text{rb}}, \gr{\phi}_1^{\text{rb}}, \hdots, \gr{\phi}_N^{\text{rb}}\},
\qquad
Q_{\rb} = \text{span}\{\psi_1^{\text{rb}}, \hdots, \psi_N^{\text{rb}}\},
\end{equation*}
for some integer $N < N_T$, and $\gr{g}^{\rb} \equiv \langle\gr{u}_h\rangle$. The FE degrees of freedom of the bases in $\gr{V}_{\rb}$ and $Q_{\rb}$ are also stored as columns of the basis functions matrices $Z_{\rb}^{\gr{u}}$ and $Z_{\rb}^{p}$ for velocity and pressure, respectively.

\subsection{Evaluation of the ROM by Galerkin projection}\label{subsec:Gproj}
During the evaluation of the ROM, one seeks a solution of the form
\begin{equation*}
\gr{u}_{\rb}(t)= \gr{g}^{\rb} + \sum_{n=1}^{N}\left(u_{\rb}(t)\right)_{n}\gr{\varphi}_n^{\rb} + \sum_{n=1}^{N}\left(u_{\rb}(t)\right)_{n + N}\gr{\phi}_n^{\rb},
\qquad
p_{\rb}(t)= \sum_{n=1}^{N}\left(p_{\rb}(t)\right)_{n}\psi_n^{\rb}
\end{equation*}
through a Galerkin projection. For the sake of notation we define $\mathbf{u}_{\rb}\left(t\right)=\left(1, \left(u_{\rb}(t)\right)_{1},\dots,\left(u_{\rb}(t)\right)_{2N}\right)^T \in \mathbb{R}^{2N+1}$ and $\mathbf{p}_{\rb}\left(t\right)=\left(\left(p_{\rb}(t)\right)_{1},\dots,\left(p_{\rb}(t)\right)_{N}\right)^T \in \mathbb{R}^N$ the vector of coefficients of the reduced solution.

Following \cite{Shafqat,PACCIARINI20141} we present here two approaches for the evaluation of the ROM, that we will refer here as \emph{consistent ROM} and \emph{non-consistent ROM}. In the consistent ROM, a Galerkin projection of the VMS Navier-Stokes equations \eqref{NS coarse} onto the reduced basis spaces is carried out. In contrast, in the non-consistent ROM projection of the standard Navier-Stokes equations \eqref{NS weak discrete} onto the reduced basis is sought. As snapshots were generated from \eqref{NS coarse}, the latter choice features a modelling inconsistency between the full order and reduced order levels, hence the name.

\subsubsection{Consistent ROM}
As mentioned, the aim of the \emph{consistent ROM} is to provide a reduced order model which is fully consistent with the full order model, thus including the projection of all stabilization terms of the VMS method. The resulting ROM requires therefore the solution of the following dynamical system:\footnote{In a practical implementation the row/column of the system corresponding to the first entry of $\mathbf{u}_{\rb}$ (which is a known coefficient that accounts for the non-homogeneous boundary condition) shall be condensed.}
\begin{equation*}
\begin{cases}
\begin{split}
M_{\rb}\frac{d\mathbf{u}_{\rb}(t)}{dt}&+A_{\rb}\mathbf{u}_{\rb}(t)+C_{\rb}(\mathbf{u}_{\rb}(t))\mathbf{u}_{\rb}(t)
+D_{\rb}(\mathbf{u}_{\rb}(t),\mathbf{p}_{\rb}(t))
+B^T_{\rb}\mathbf{p}_{\rb}(t)
=\gr{F}_{\rb}(t),\\
B_{\rb}\mathbf{u}_{\rb}(t)&=E_{\rb}(\mathbf{u}_{\rb}(t),\mathbf{p}_{\rb}(t)),
\end{split}
\end{cases}
\end{equation*}
being
\begin{align*}
&
M_{\rb} = {Z_{\rb}^{\gr{u}}}^T M Z_{\rb}^{\gr{u}},\qquad
A_{\rb} = {Z_{\rb}^{\gr{u}}}^T A Z_{\rb}^{\gr{u}},\qquad
B_{\rb} = {Z_{\rb}^{p}}^T B Z_{\rb}^{\gr{u}},\\
&
C_{\rb}\left(\mathbf{u}_{\rb}\right) = {Z_{\rb}^{\gr{u}}}^T C\left(Z_{\rb}^{\gr{u}} \mathbf{u}_{\rb}(t)\right) Z_{\rb}^{\gr{u}},\qquad
D_{\rb}(\mathbf{u}_{\rb}(t),\mathbf{p}_{\rb}(t)) = {Z_{\rb}^{\gr{u}}}^T D(Z_{\rb}^{\gr{u}} \mathbf{u}_{\rb}(t),Z_{\rb}^{p}\mathbf{p}_{\rb}(t)),\\
&
E_{\rb}(\mathbf{u}_{\rb}(t),\mathbf{p}_{\rb}(t)) = {Z_{\rb}^{p}}^T E(Z_{\rb}^{\gr{u}} \mathbf{u}_{\rb}(t),Z_{\rb}^{p}\mathbf{p}_{\rb}(t)),\qquad
F_{\rb} = {Z_{\rb}^{\gr{u}}}^T F.
\end{align*}
As the terms $D_{\rb}(\mathbf{u}_{\rb}(t),\mathbf{p}_{\rb}(t))$ and $E_{\rb}(\mathbf{u}_{\rb}(t),\mathbf{p}_{\rb}(t))$, arising from VMS stabilization {(as defined in \autoref{subsec:FEM_VMS})}, are highly nonlinear, their evaluation should be handled with care to preserve the efficiency of the resulting ROM, enforcing the so-called offline-online decoupling (see e.g. \cite{libroRozza}). Hyper reduction techniques, such as empirical interpolation \cite{BARRAULT2004667}, GNAT \cite{Carlberg2013623} or Gappy-POD \cite{Willcox2006} methods are available to this end. Since the goal of this preliminary investigation is to compare the accuracy and feasibility of the two proposed ROMs, without aiming at the greatest amount of efficiency, we will allow an inefficient evaluation of the nonlinear terms. {In order to ensure consistency, the same values for the stabilization coefficients appearing in \eqref{fine scale explicit} in the FOM are used for the ROM. Alternative options are explored in literature for the modeling of such coefficients at the reduced order level \cite{Giere2015}, and would result in an intermediate choice between our two proposed methods. Such choice would result in a ROM with stabilization, yet not full consistency with the FOM. A detailed comparison of all the resulting options is out of the scope of the current work, and will be pursued in future.}

\subsubsection{Non-consistent ROM}
In contrast, the goal of the \emph{non-consistent ROM} is to perform a Galerkin projection of the standard Navier-Stokes equations without any additional stabilization term. The resulting dynamical system is:
\begin{equation*}
\begin{cases}
\begin{split}
M_{\rb}\frac{d\mathbf{u}_{\rb}(t)}{dt}&+A_{\rb}\mathbf{u}_{\rb}(t)+C_{\rb}(\mathbf{u}_{\rb}(t))\mathbf{u}_{\rb}(t)
+B^T_{\rb}\mathbf{p}_{\rb}(t)
=\gr{F}_{\rb}(t),\\
B_{\rb}\mathbf{u}_{\rb}(t)&=\mathbf{0}.
\end{split}
\end{cases}
\end{equation*}
This ROM seems attractive since it only has a nonlinear contribution relative to the convective term, which is a quadratic nonlinearity.
In this case an efficient offline-online decoupling may be easily obtained with the aid of a third tensor precomputed at the end of the construction stage \cite{Ballarin2015}. This model is however not consistent with the full order model used to generate the snapshots because the VMS stabilization terms are neglected during the projection stage. Yet, such inconsistent options are being used with success in literature, as different models at the full order and reduced order levels are sought in \cite{Bergmann2009516,wang_turb}.

\section{A numerical example: flow around a circular cylinder}\label{sec:numer_example}

\begin{figure}
\begin{minipage}{\textwidth}
\begin{minipage}{0.33\textwidth}
\includegraphics[width=\textwidth]{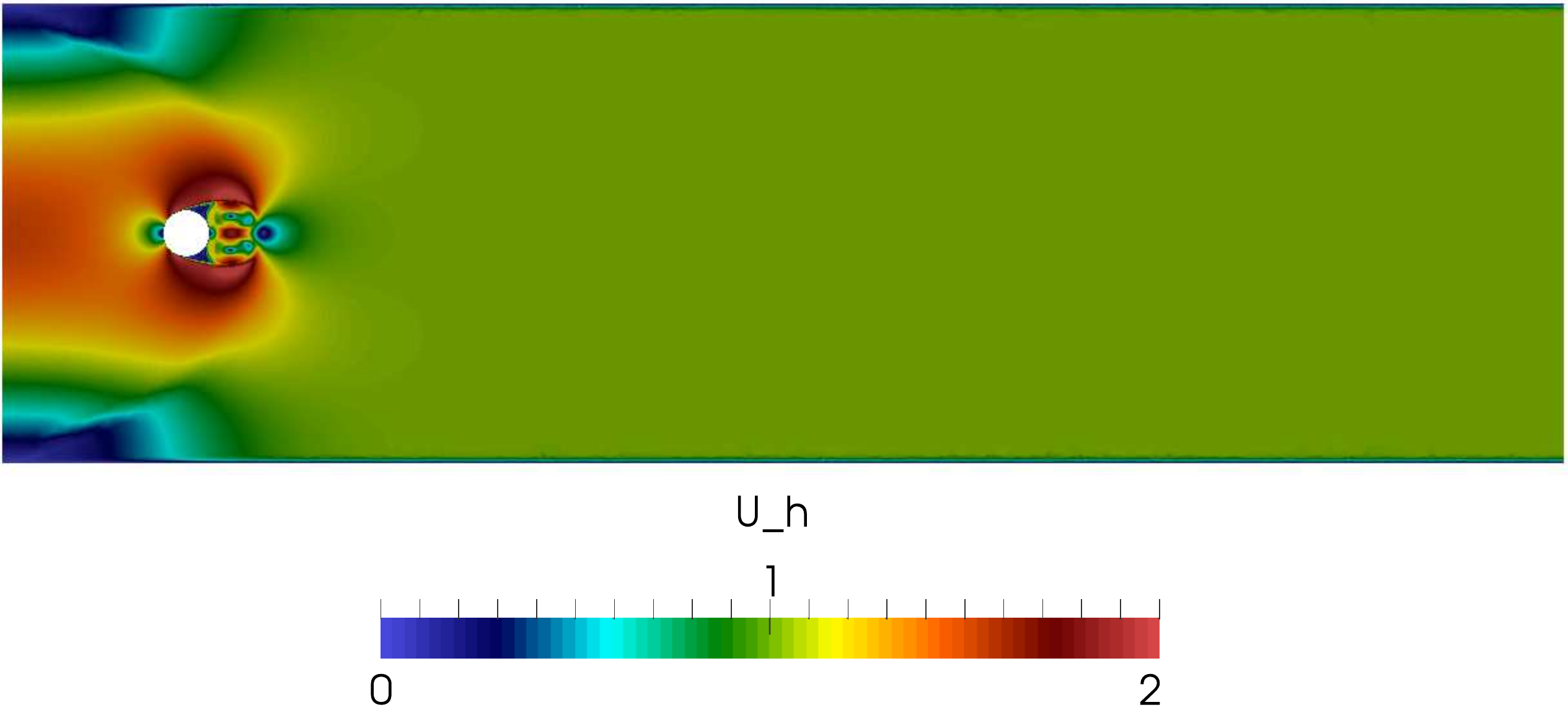}
\end{minipage}
\begin{minipage}{0.33\textwidth}
\includegraphics[width=\textwidth]{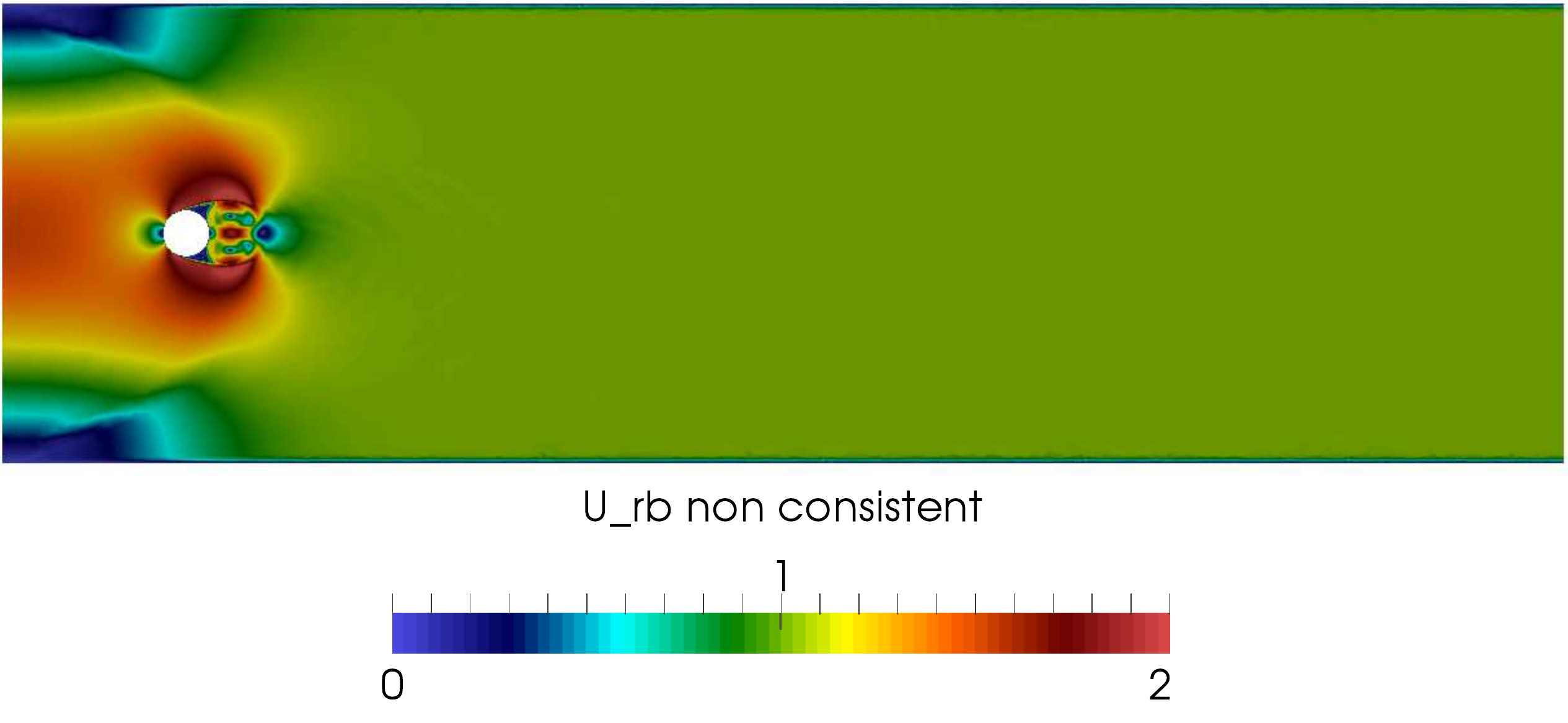}
\end{minipage}
\begin{minipage}{0.33\textwidth}
\includegraphics[width=\textwidth]{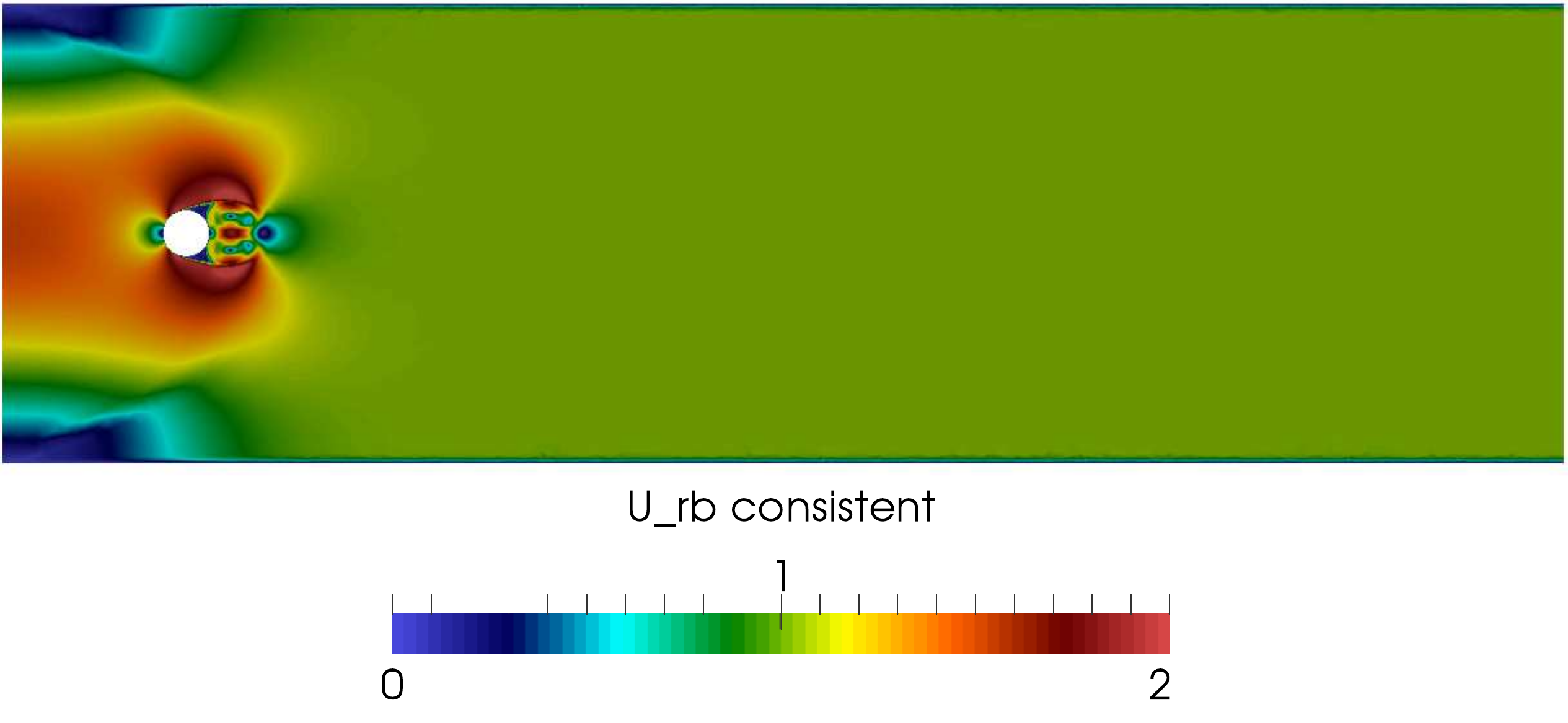}
\end{minipage}
\begin{minipage}{0.33\textwidth}
\includegraphics[width=\textwidth]{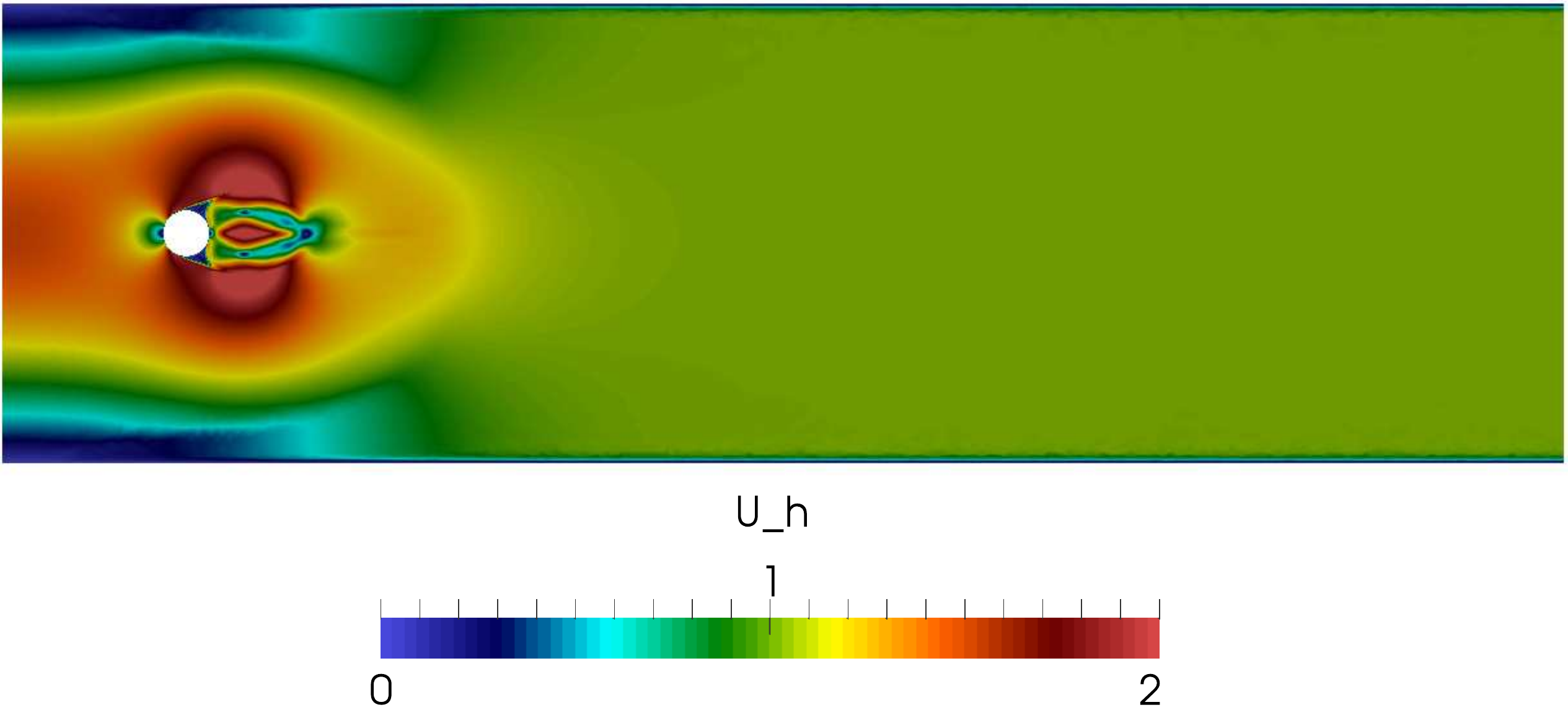}
\end{minipage}
\begin{minipage}{0.33\textwidth}
\includegraphics[width=\textwidth]{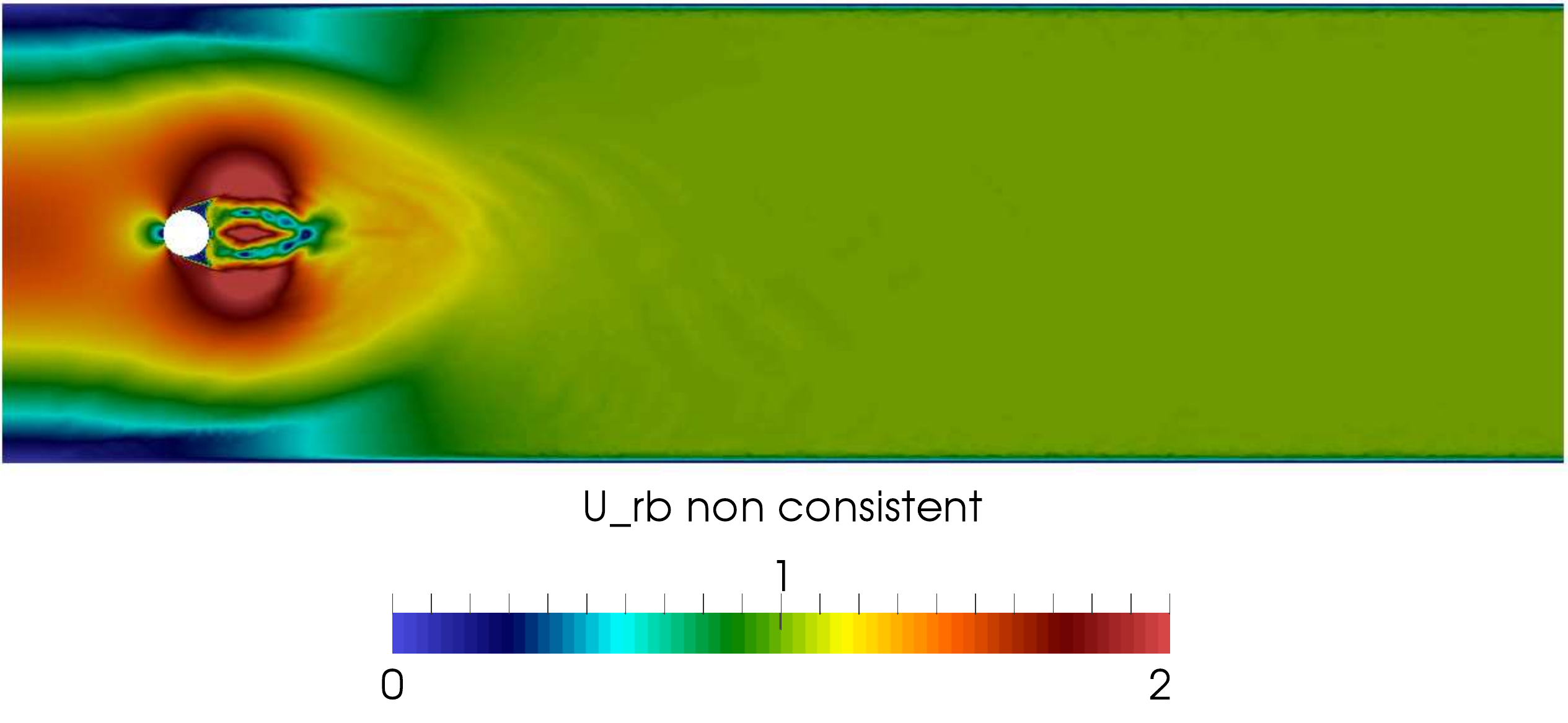}
\end{minipage}
\begin{minipage}{0.33\textwidth}
\includegraphics[width=\textwidth]{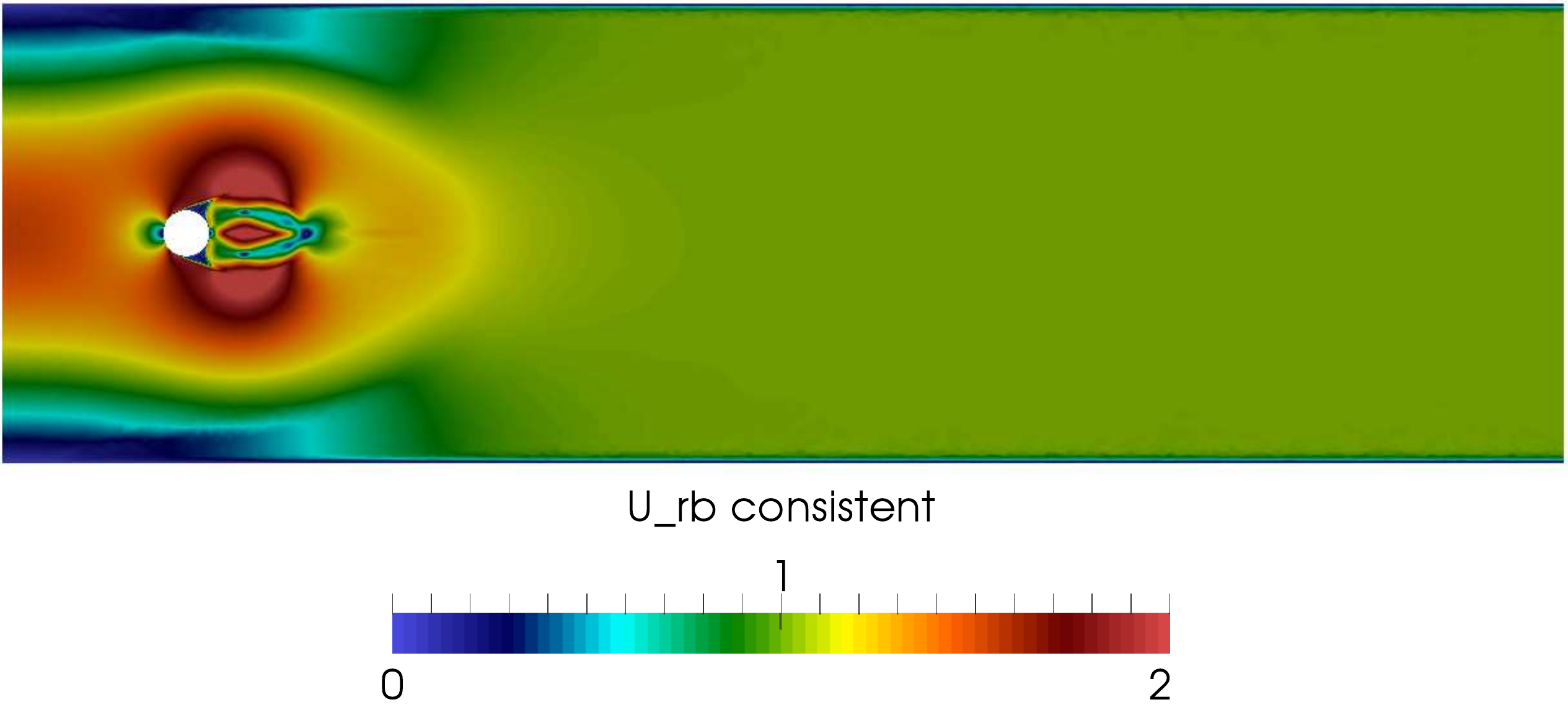}
\end{minipage}
\begin{minipage}{0.33\textwidth}
\includegraphics[width=\textwidth]{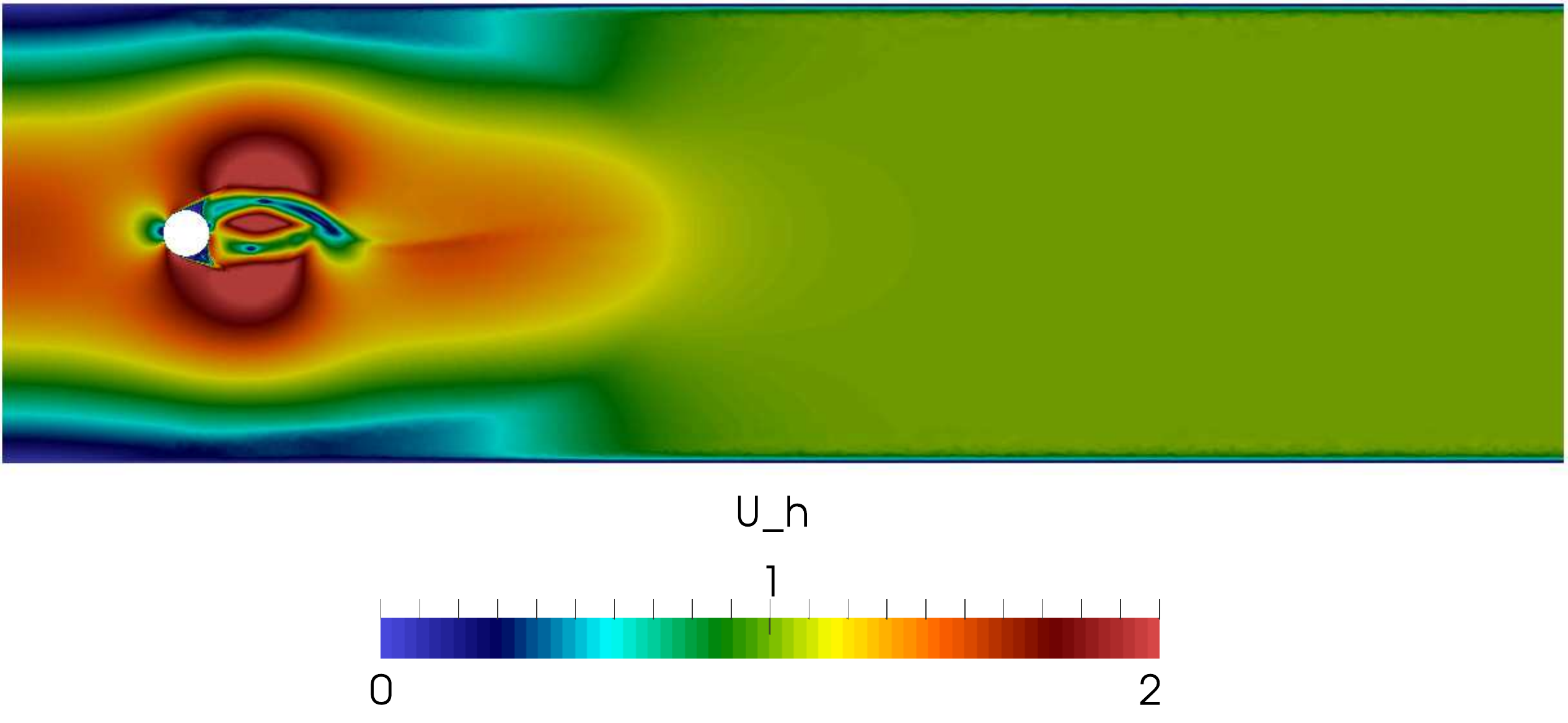}
\end{minipage}
\begin{minipage}{0.33\textwidth}
\includegraphics[width=\textwidth]{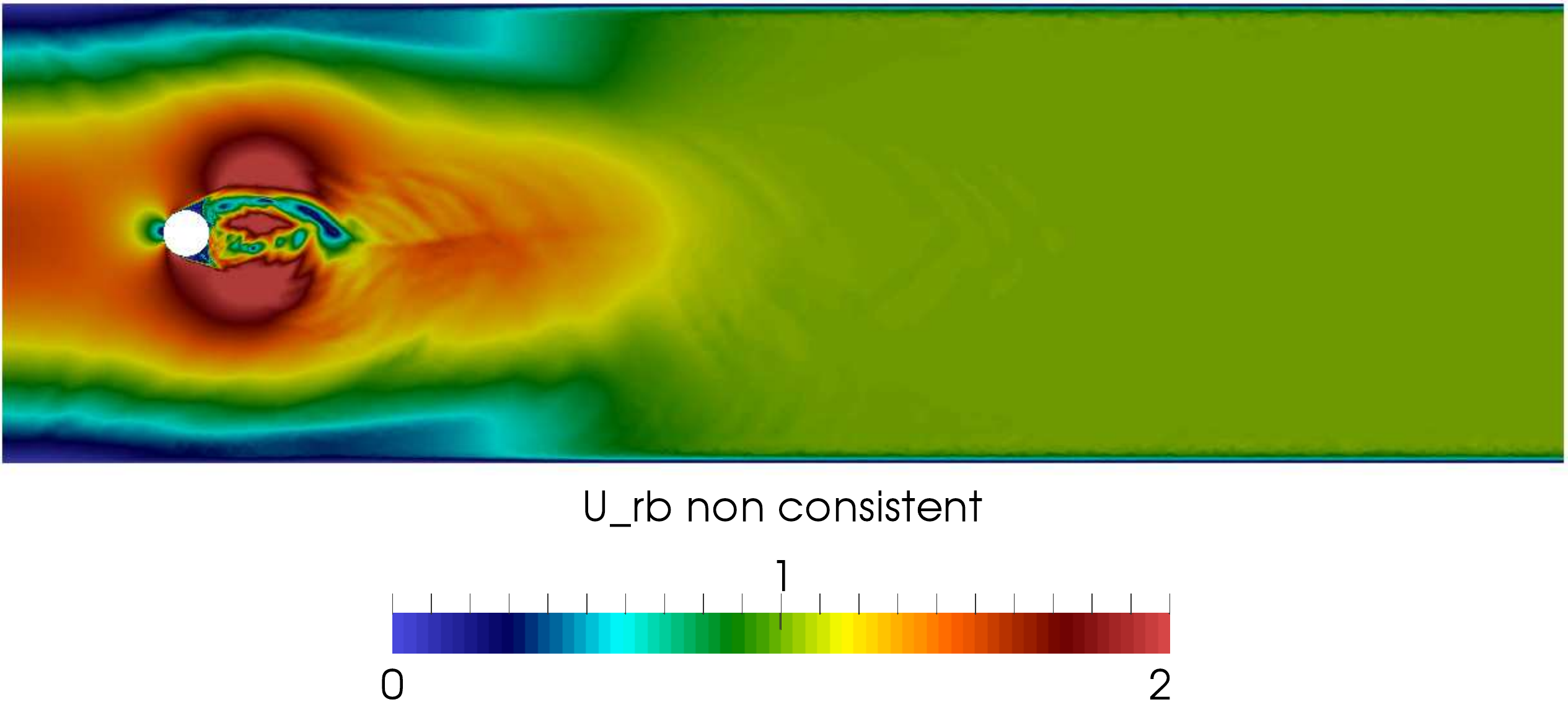}
\end{minipage}
\begin{minipage}{0.33\textwidth}
\includegraphics[width=\textwidth]{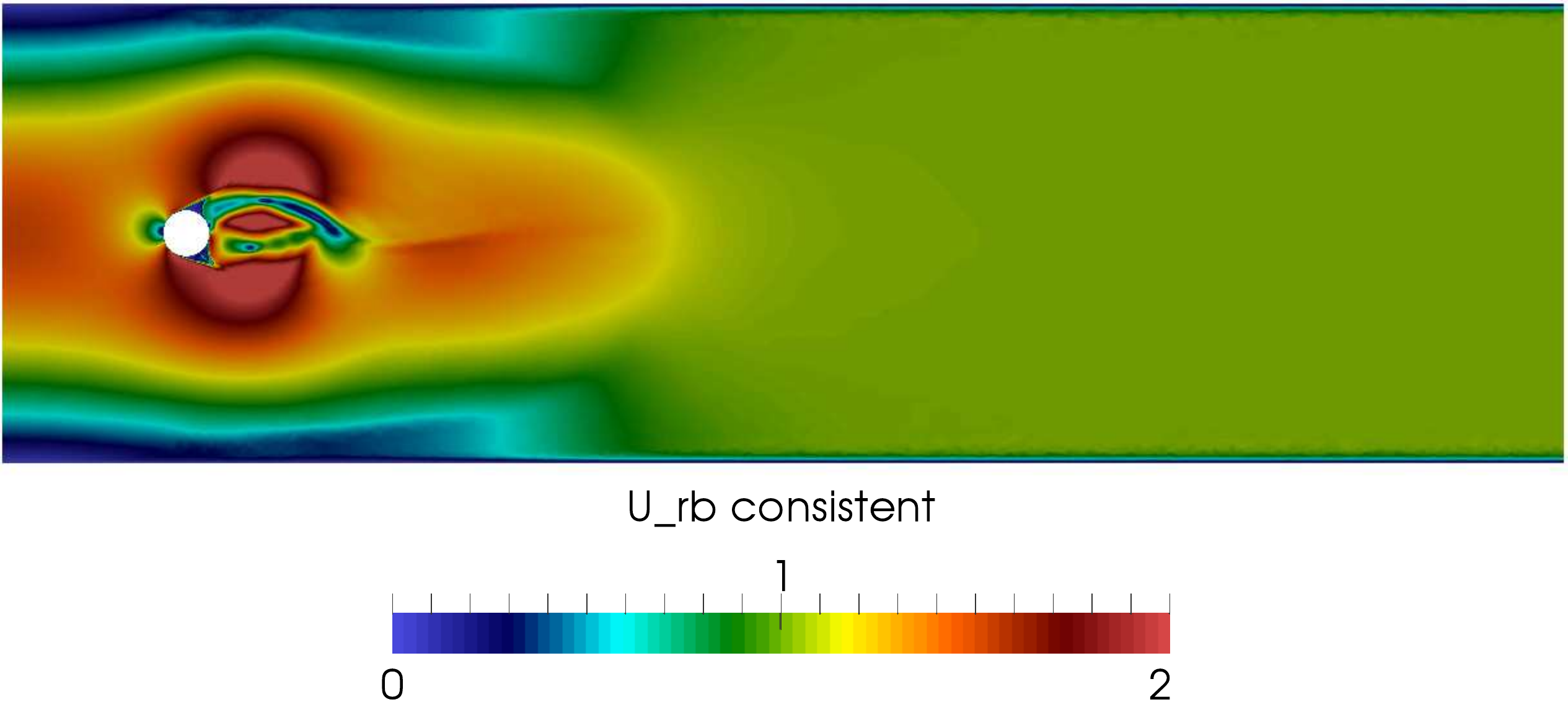}
\end{minipage}
\begin{minipage}{0.33\textwidth}
\includegraphics[width=\textwidth]{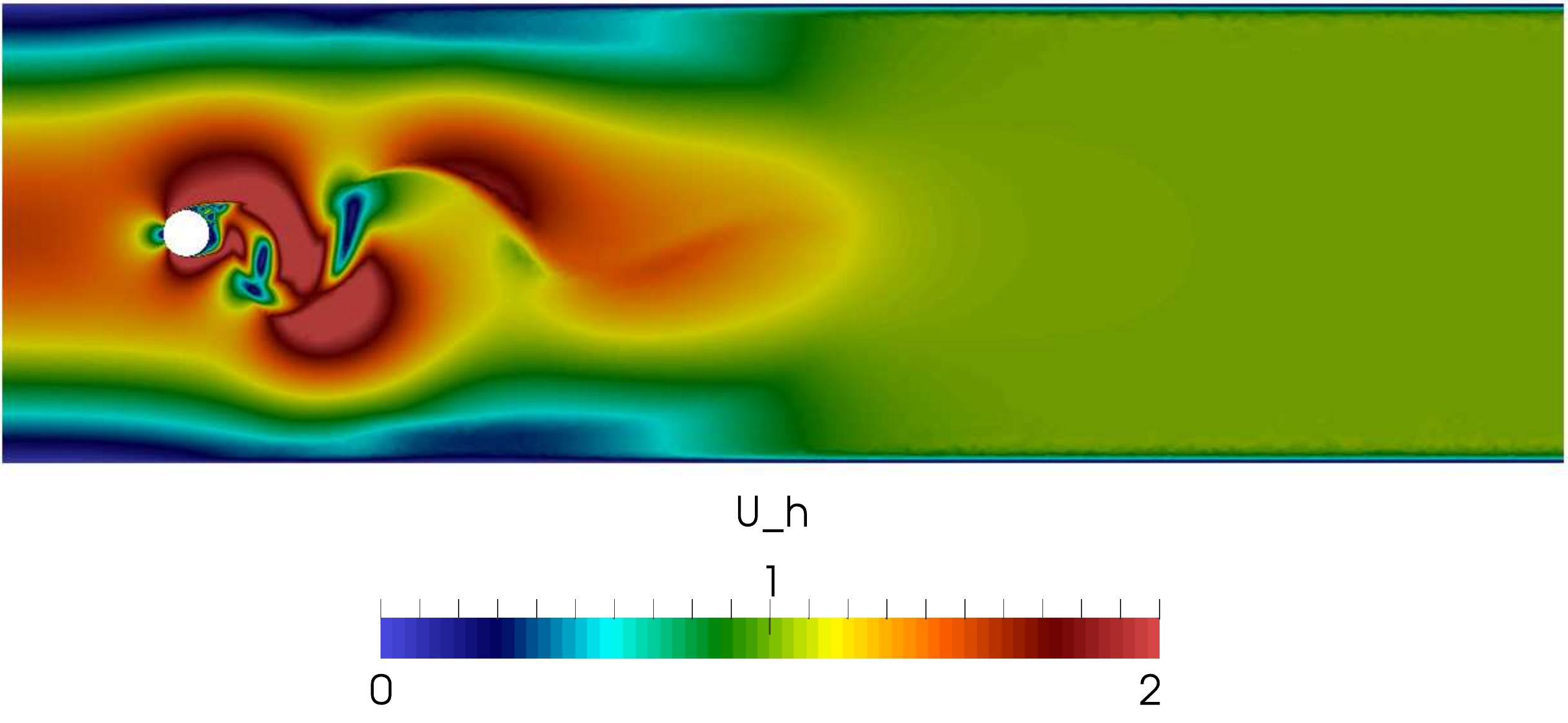}
\end{minipage}
\begin{minipage}{0.33\textwidth}
\includegraphics[width=\textwidth]{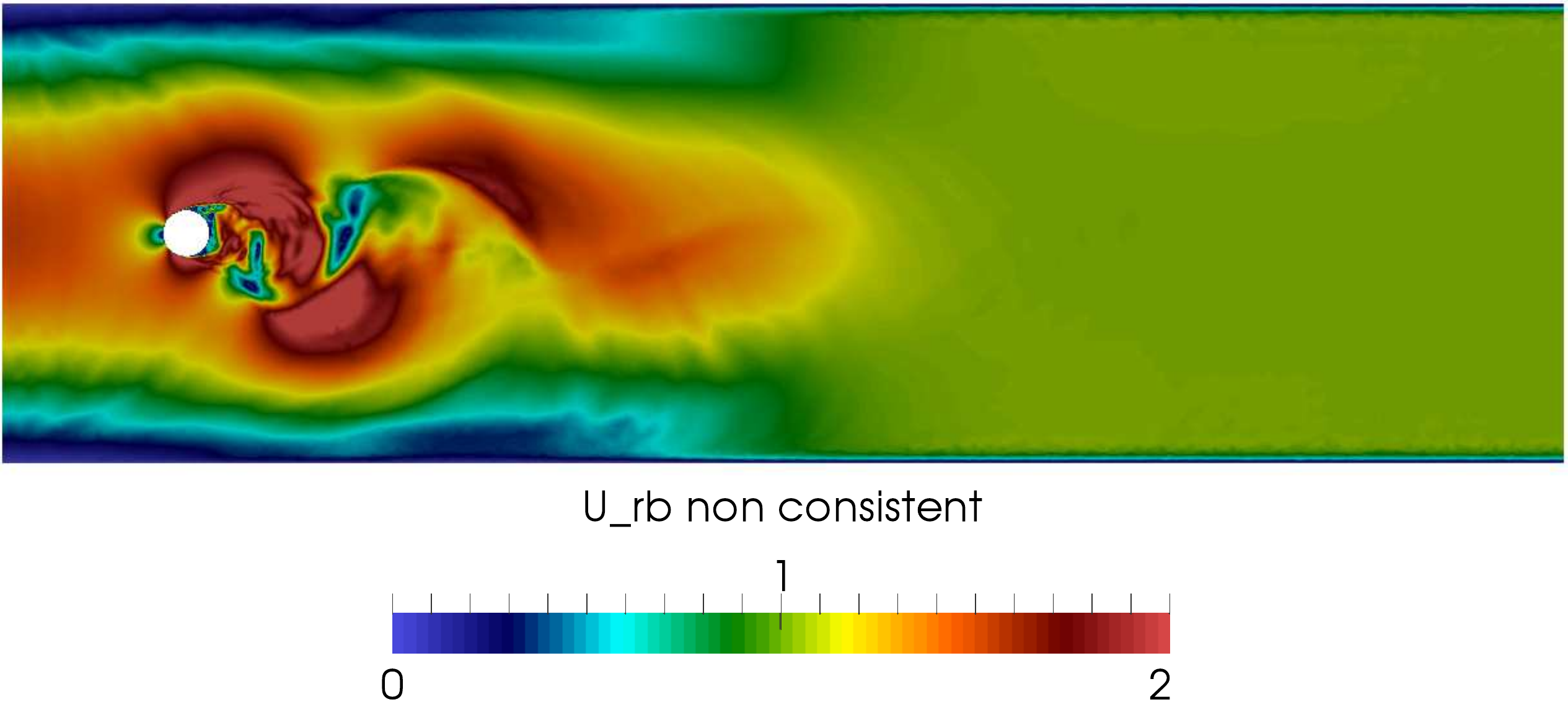}
\end{minipage}
\begin{minipage}{0.33\textwidth}
\includegraphics[width=\textwidth]{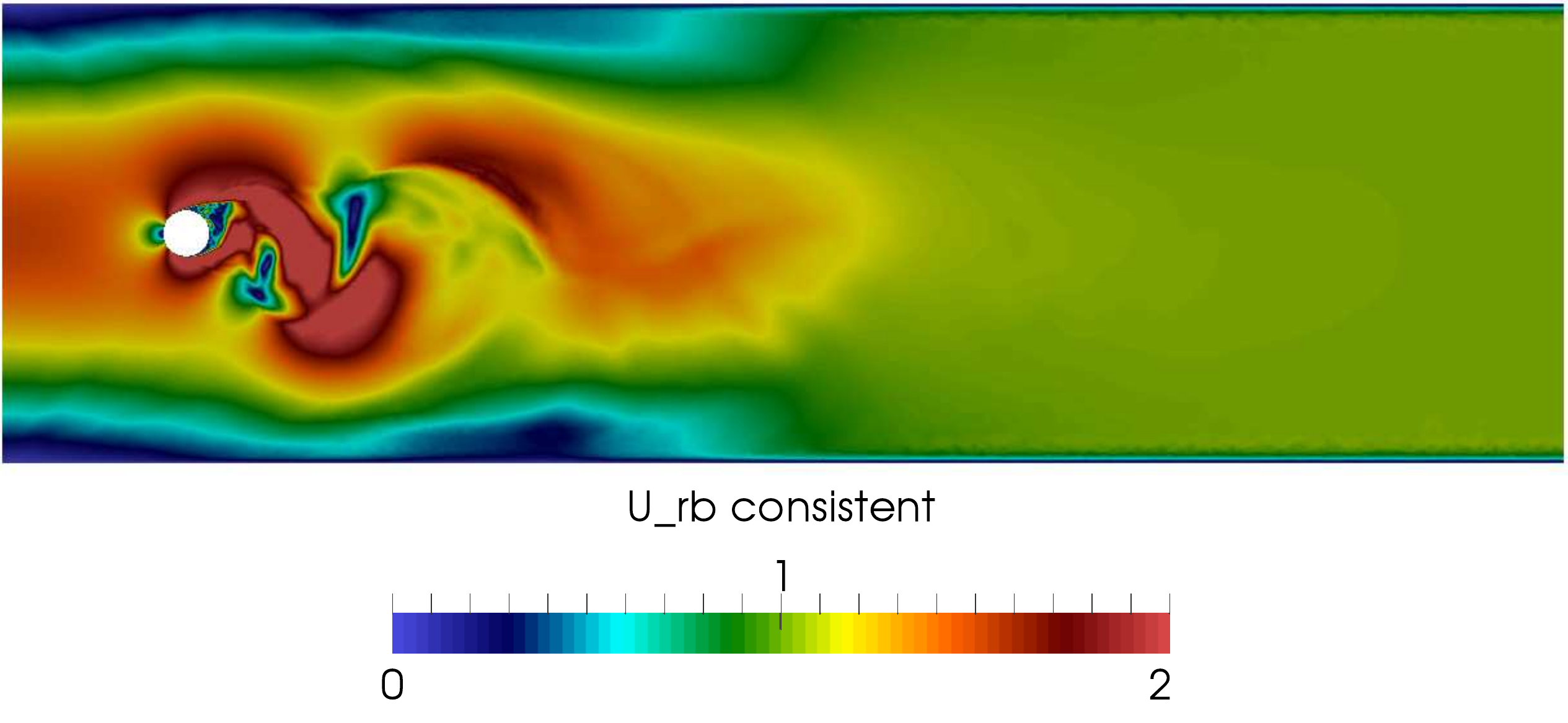}
\end{minipage}
\end{minipage}
\caption{Comparison of the flow fields in terms of velocity magnitude for $t=0.5$ s, $t=1.0$ s, $t=1.5$ s, $t=2.0$ s (from top to bottom) for full order model, non-consistent ROM with supremizer, consistent ROM with supremizer (from left to right). ROM evaluations were carried out for $N=100$.}
\label{fig:plots_u}
\end{figure}

\begin{figure}
\begin{minipage}{\textwidth}
\begin{minipage}{0.33\textwidth}
\includegraphics[width=\textwidth]{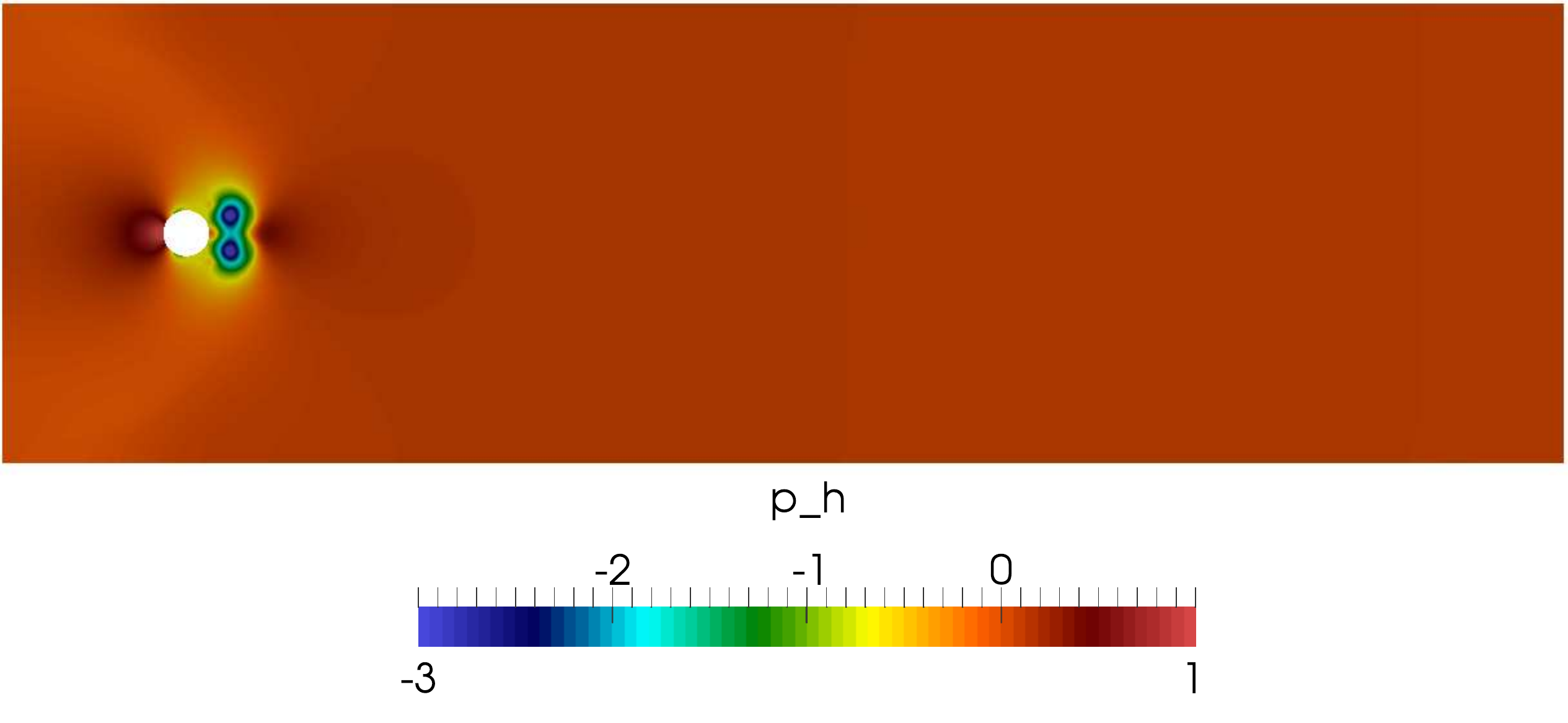}
\end{minipage}
\begin{minipage}{0.33\textwidth}
\includegraphics[width=\textwidth]{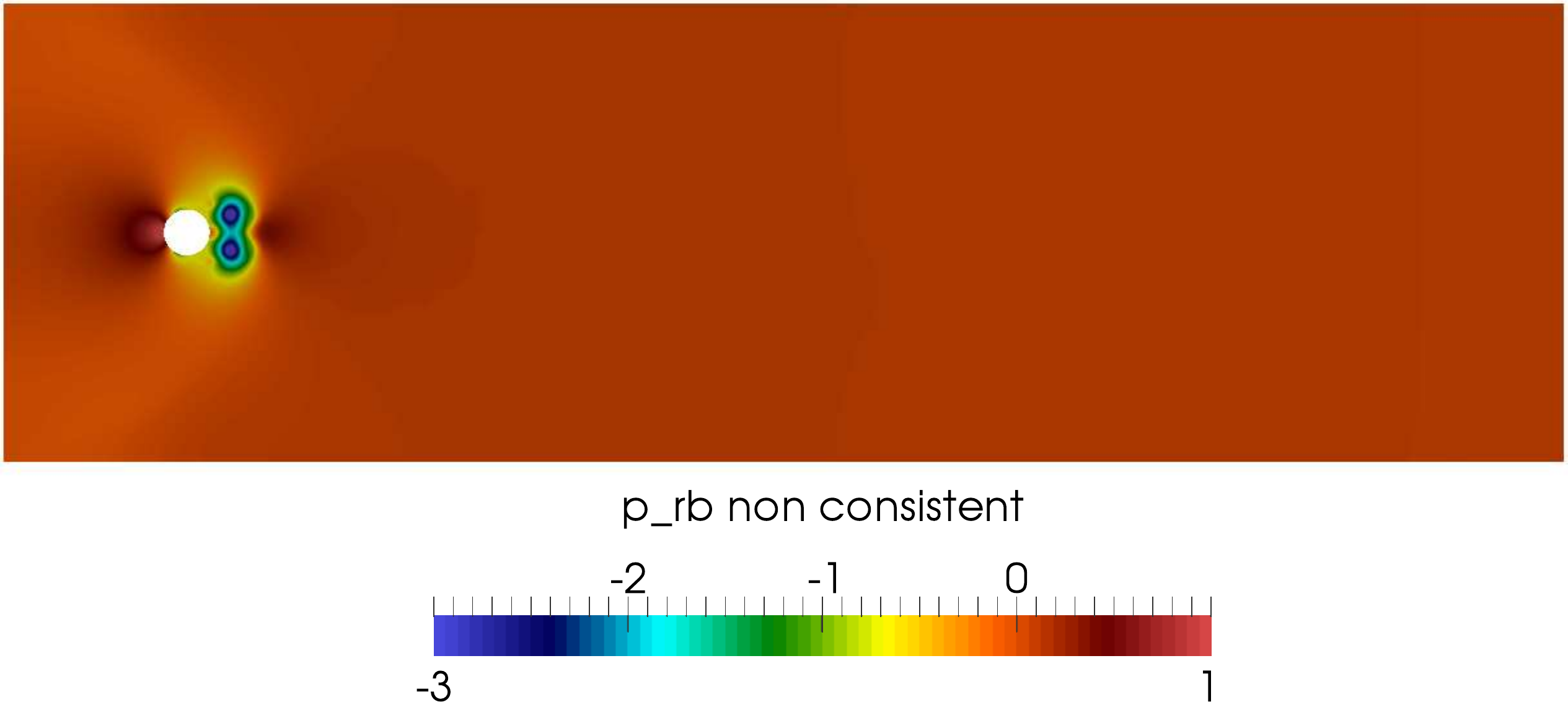}
\end{minipage}
\begin{minipage}{0.33\textwidth}
\includegraphics[width=\textwidth]{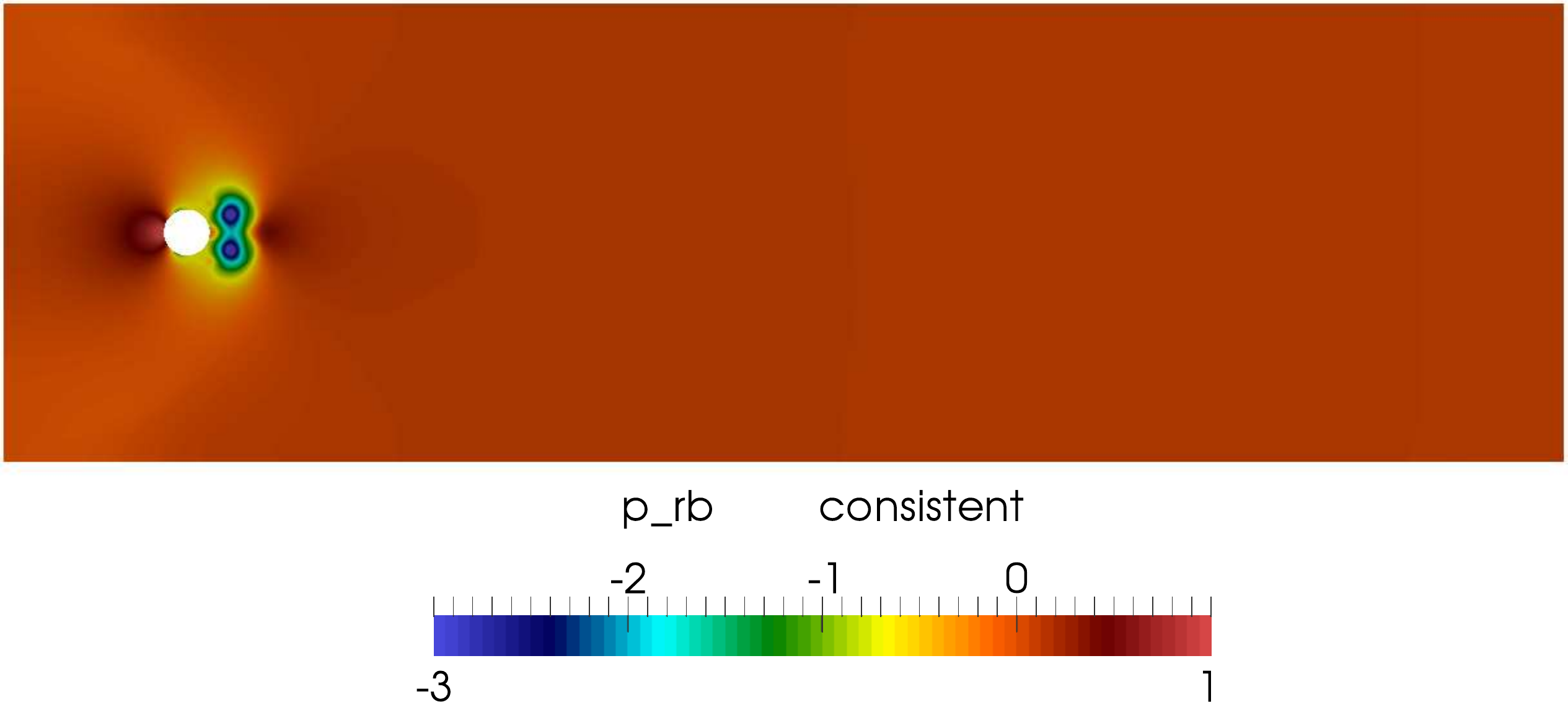}
\end{minipage}
\begin{minipage}{0.33\textwidth}
\includegraphics[width=\textwidth]{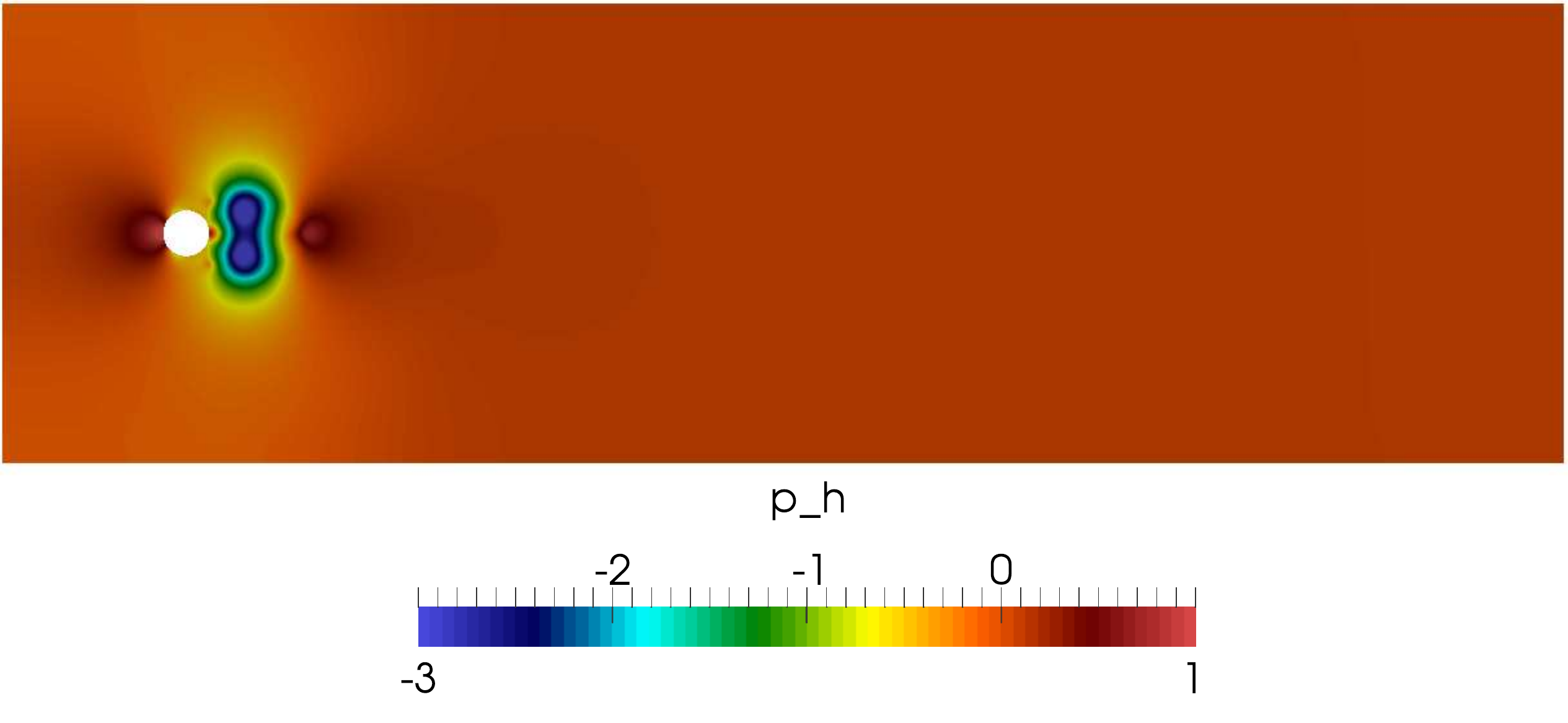}
\end{minipage}
\begin{minipage}{0.33\textwidth}
\includegraphics[width=\textwidth]{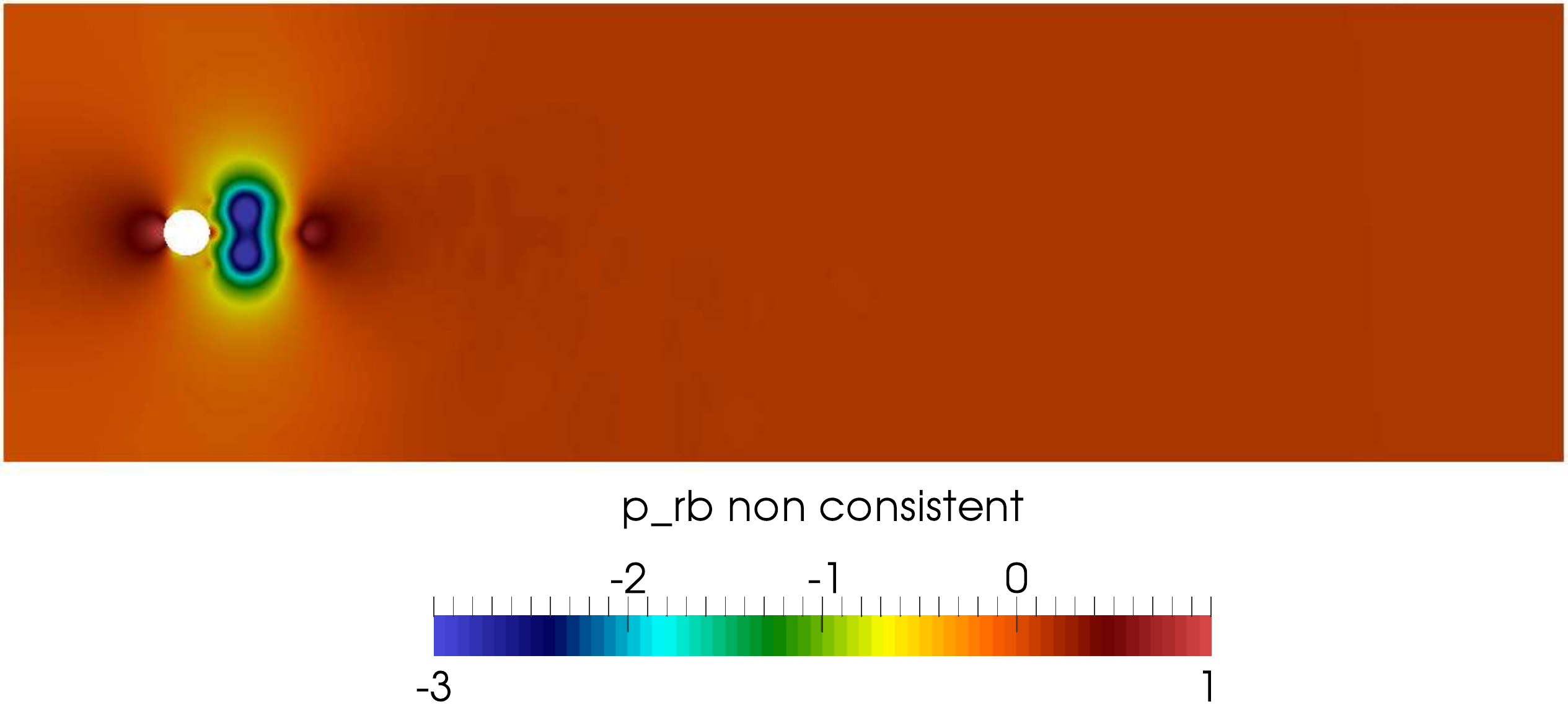}
\end{minipage}
\begin{minipage}{0.33\textwidth}
\includegraphics[width=\textwidth]{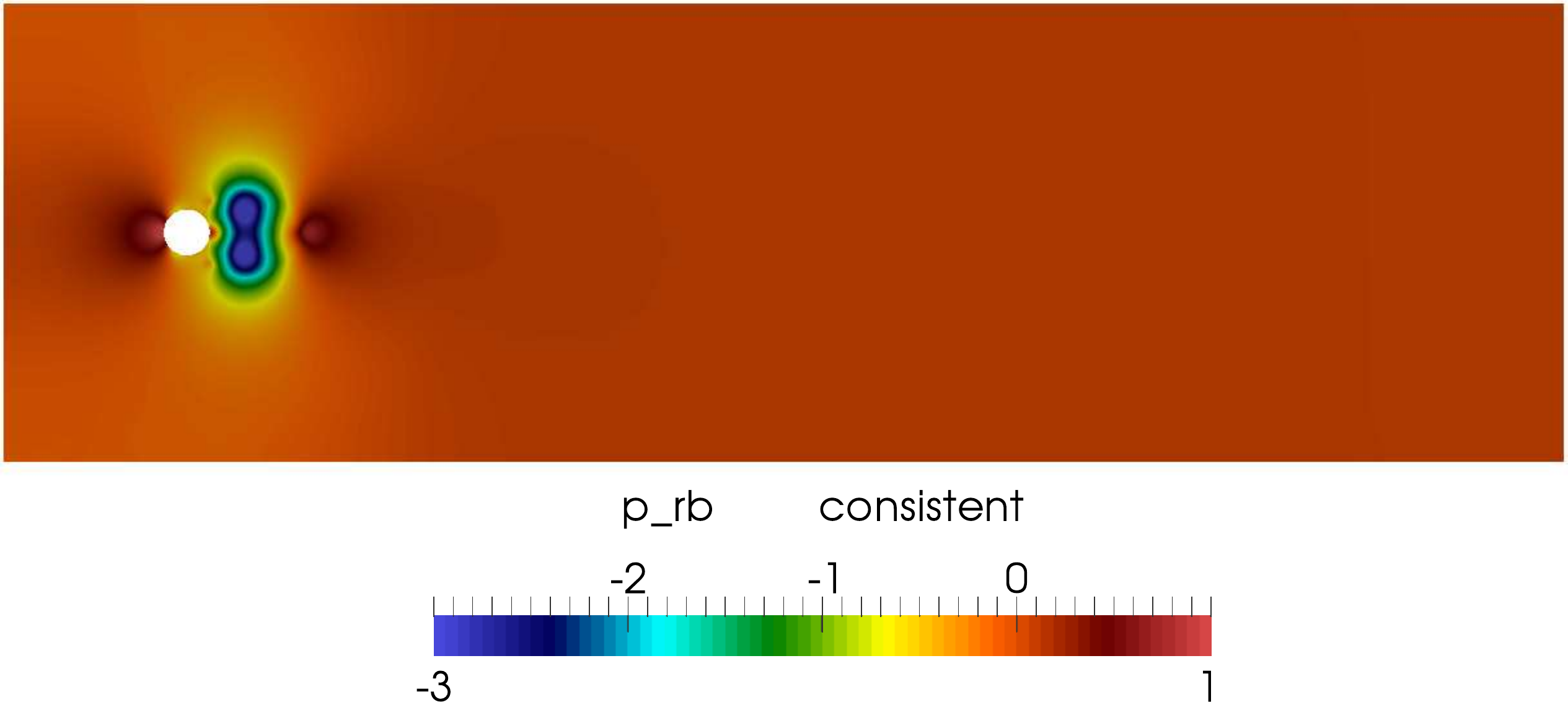}
\end{minipage}
\begin{minipage}{0.33\textwidth}
\includegraphics[width=\textwidth]{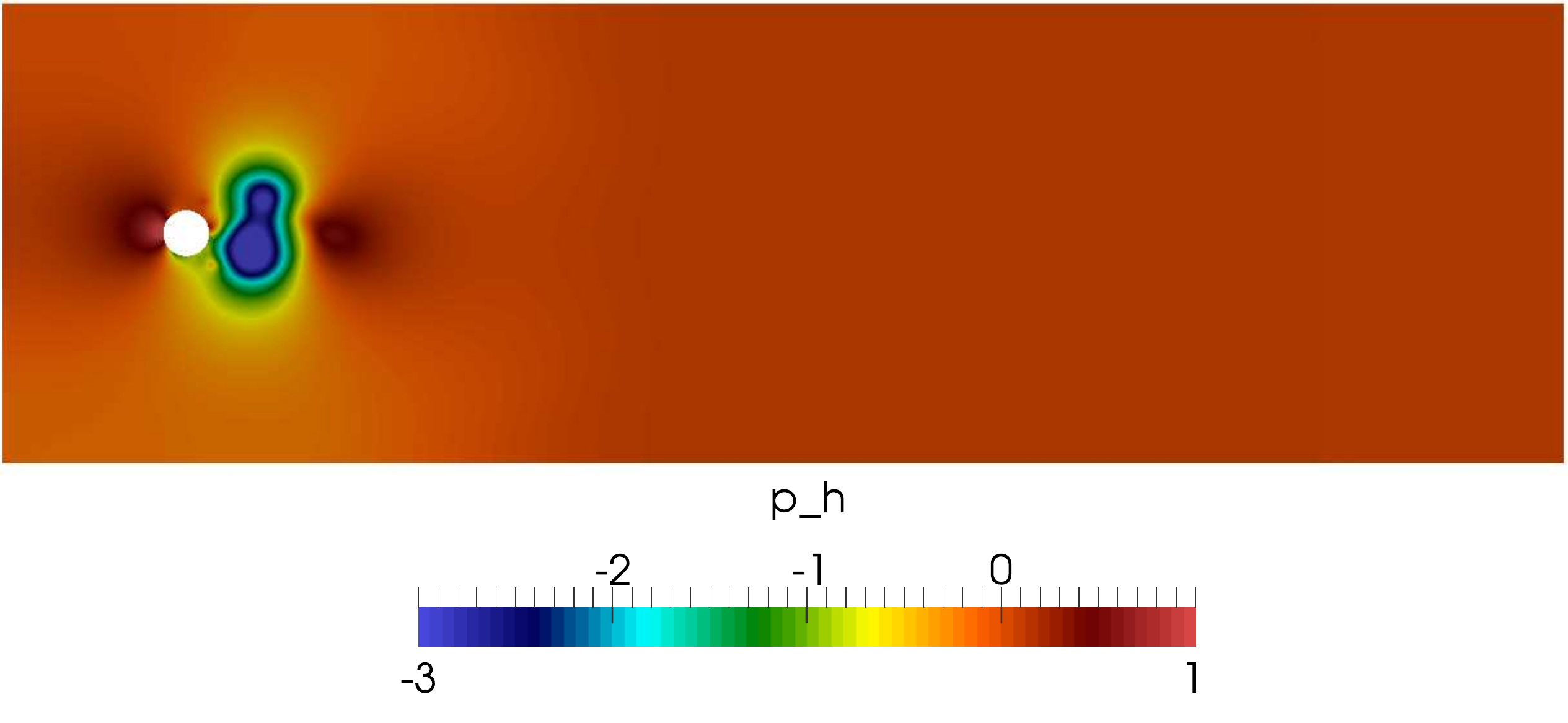}
\end{minipage}
\begin{minipage}{0.33\textwidth}
\includegraphics[width=\textwidth]{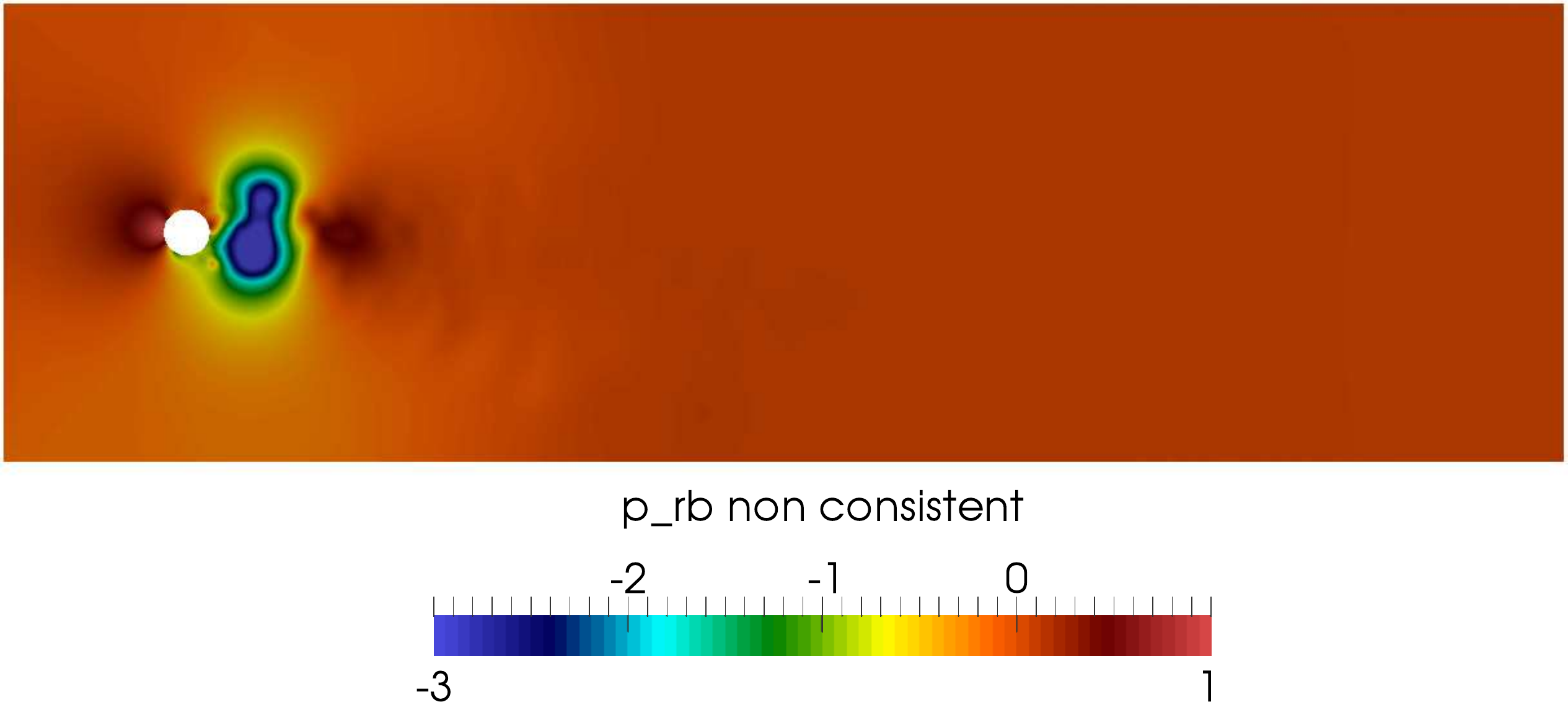}
\end{minipage}
\begin{minipage}{0.33\textwidth}
\includegraphics[width=\textwidth]{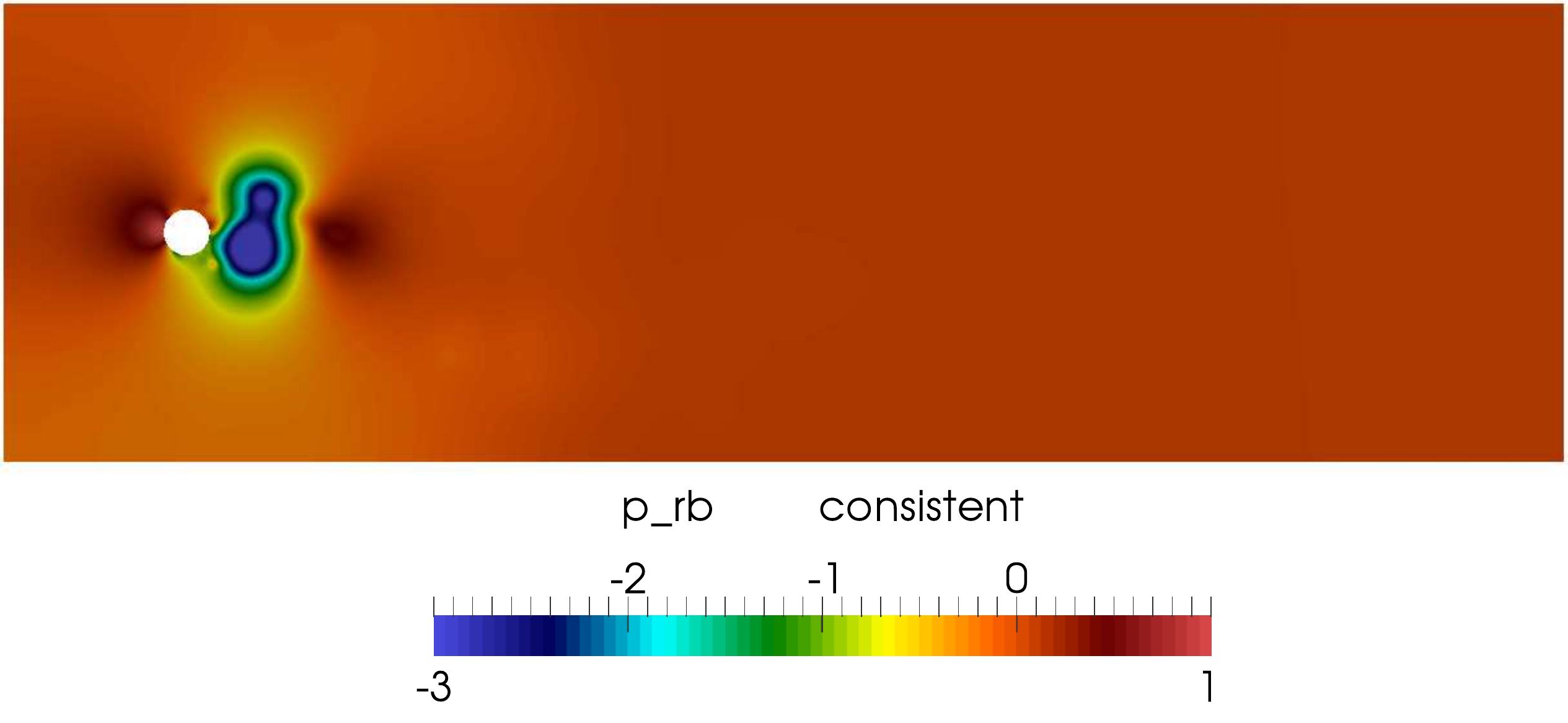}
\end{minipage}
\begin{minipage}{0.33\textwidth}
\includegraphics[width=\textwidth]{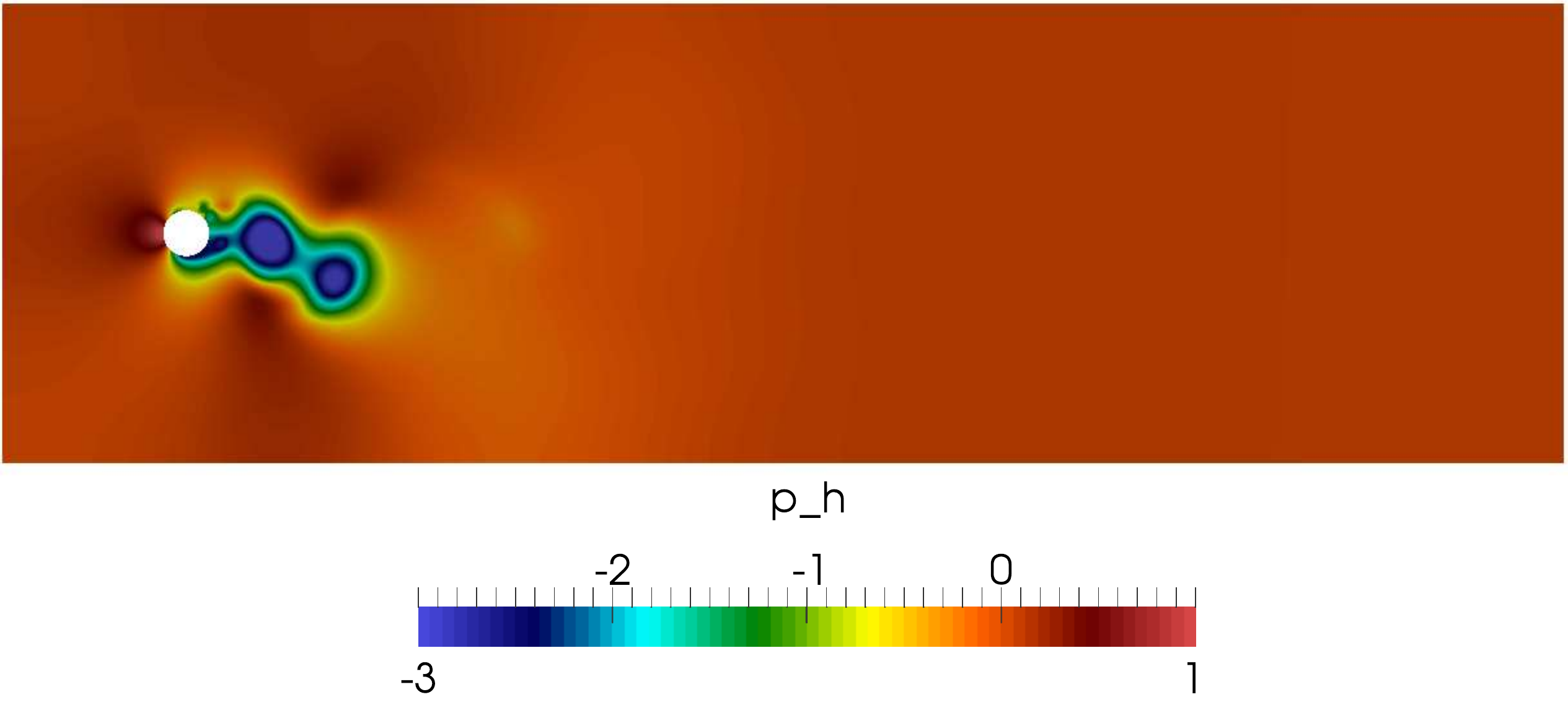}
\end{minipage}
\begin{minipage}{0.33\textwidth}
\includegraphics[width=\textwidth]{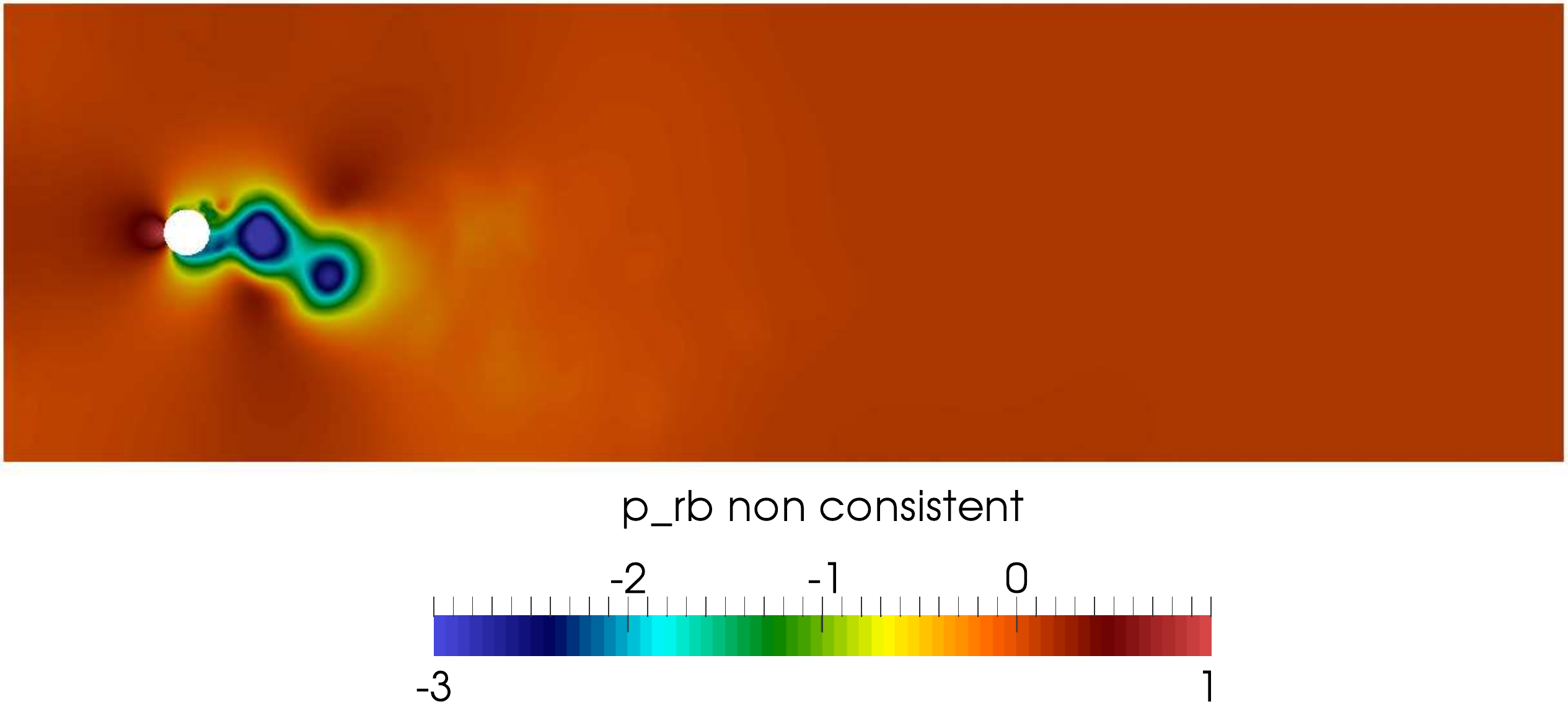}
\end{minipage}
\begin{minipage}{0.33\textwidth}
\includegraphics[width=\textwidth]{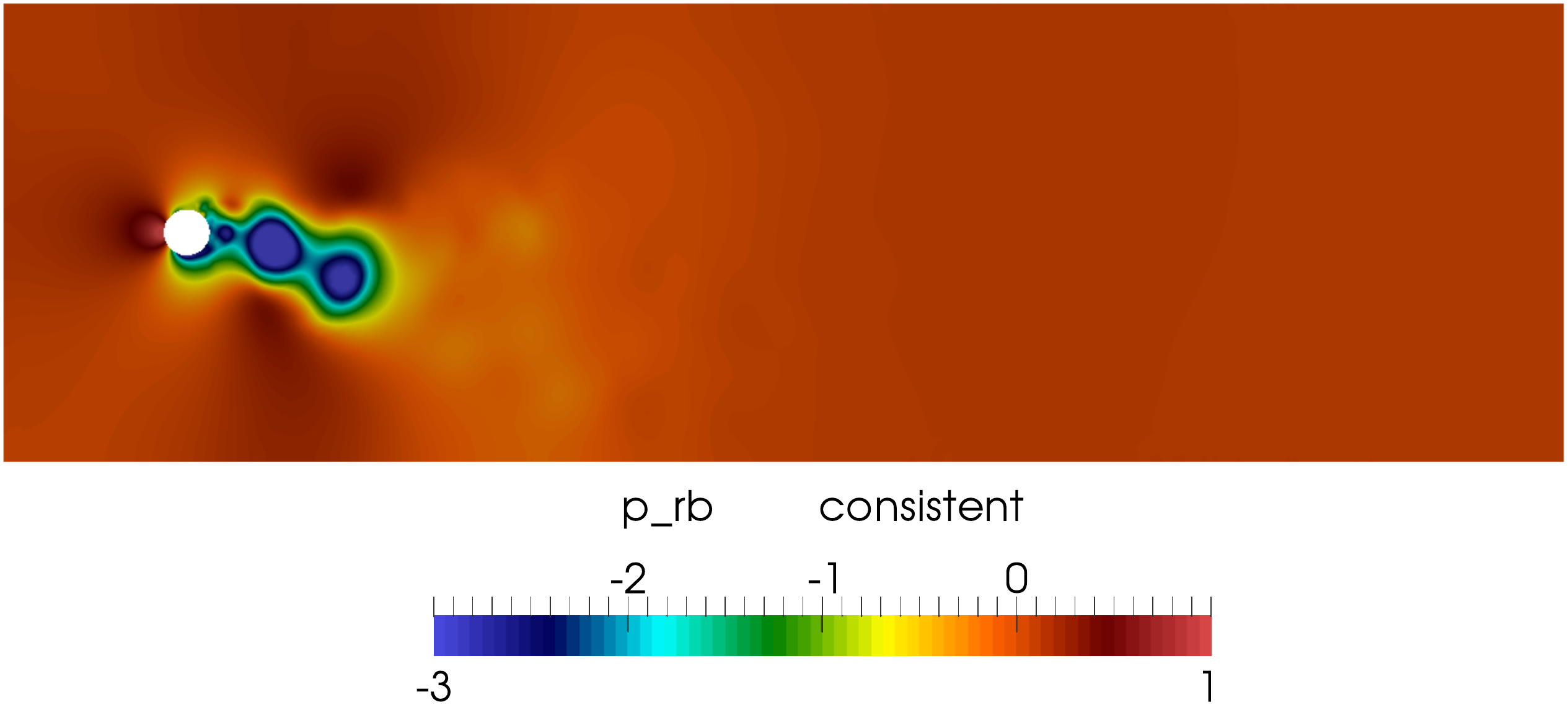}
\end{minipage}
\end{minipage}
\caption{Comparison of the flow fields in terms of pressure magnitude for $t=0.5$ s, $t=1.0$ s, $t=1.5$ s, $t=2.0$ s (from top to bottom) for full order model, non-consistent ROM with supremizer, consistent ROM with supremizer (from left to right). ROM evaluations were carried out for $N=100$.}
\label{fig:plots_p}
\end{figure}

In this section we present some numerical results on a benchmark test case involving the flow around a {2D} circular {cylinder} at moderately high Reynolds numbers (Re $\approx 5000$). The resulting flow conditions are such that a standard Navier-Stokes formulation would fail (unless being provided an extremely refined mesh). {We note that, even though turbulence is an intrinsically 3D phenomenon, {this problem is representative of a case where the numerical method does not resolve all the essential features, as is the case in turbulent flows}, and is chosen for the sake of a simpler discussion and visualization of the results (e.g. in forthcoming Figures \ref{fig:plots_u}-\ref{fig:plots_p}). Nonetheless, forthcoming papers should further test the proposed methodology on more complex 3D configurations.}

\begin{figure}
\begin{minipage}{\textwidth}
\centering
\includegraphics[width=0.8\textwidth]{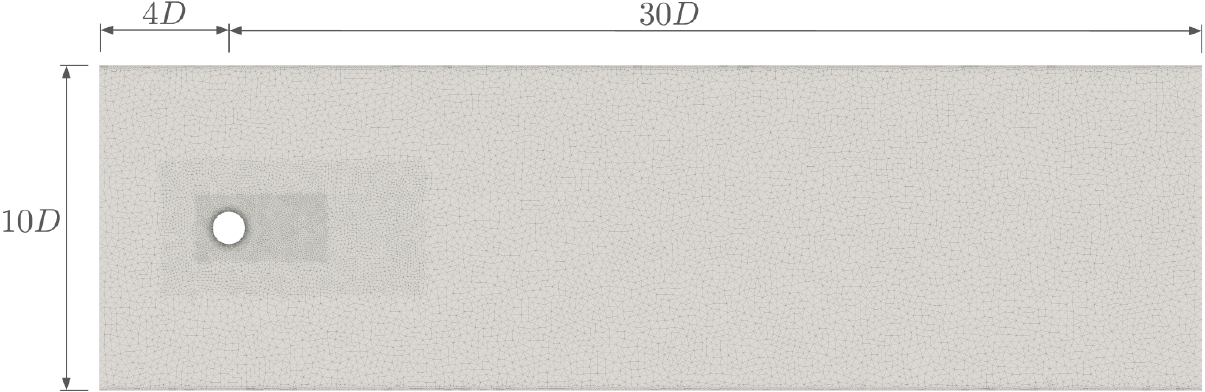}
\end{minipage}
\caption{The computational domain together with the geometrical dimensions. The diameter of the cylinder is equal to $D = 0.2\mbox{m}$}.
\label{fig:comp_mesh}
\end{figure}

The computational domain is given by a rectangular box with a circular hole. A parabolic inflow profile $\gr{g}$ is provided at the inlet, gravity is neglected $\gr{f} = \gr{0}$, and an homogeneous initial velocity is imposed. The final time is assumed to be $T = 2$ s.
The discretization of the domain is as in \autoref{fig:comp_mesh}. 
The triangular tessellation counts $37572$ cells. At the full order level, the velocity and pressure fields are approximated with $P2/P1$ finite elements, resulting in $155756$ and $20153$ degrees of freedom for velocity and pressure, respectively. A backward Euler time stepping scheme is employed, with $\Delta t = 0.002$.

\begin{figure}
\begin{minipage}{\textwidth}
\centering
\includegraphics[width=0.5\textwidth]{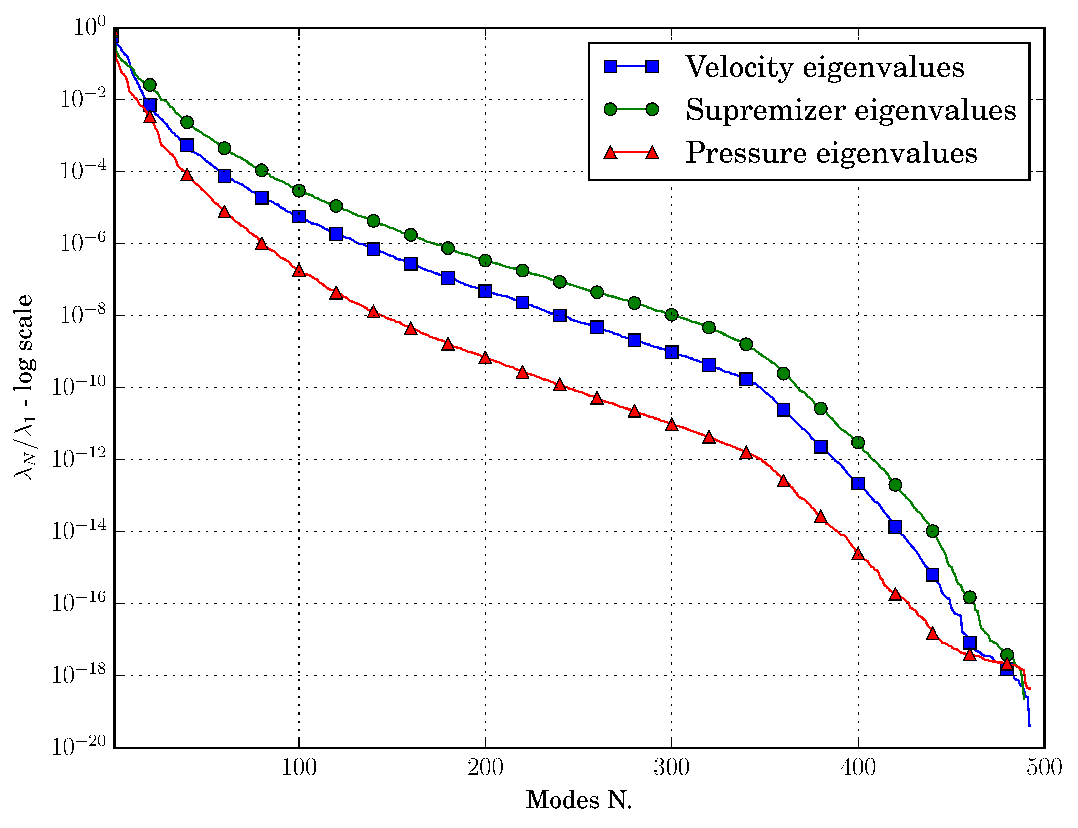}
\end{minipage}
\caption{POD eigenvalues decay for velocity, supremizer and pressure snapshots.}
\label{fig:eig_decay}
\end{figure}

The full order system \eqref{NS coarse} is solved with VMS stabilization. Undersampling is performed on the flow fields, in order to contain the size of the correlation matrices $\Sigma^{\gr{u}}$ and $\Sigma^p$, still preserving all relevant features of the transient phenomena. This results in 500 snapshots, which are used for the reduced basis generation through POD. Corresponding eigenvalues decay is shown in \autoref{fig:eig_decay} for velocity, supremizers and pressure snapshots. The eigenvalue plot, respect to a lower Reynolds number case on the same physical problem, shows a slower decay of the eigenvalue and therefore a larger number of basis functions need to be considered. This fact is due to the more complex structure of the vortices that characterizes flows with higher Reynolds number.  

\begin{figure}
\begin{minipage}{\textwidth}
\begin{minipage}{0.5\textwidth}
\centering
\includegraphics[width=\textwidth]{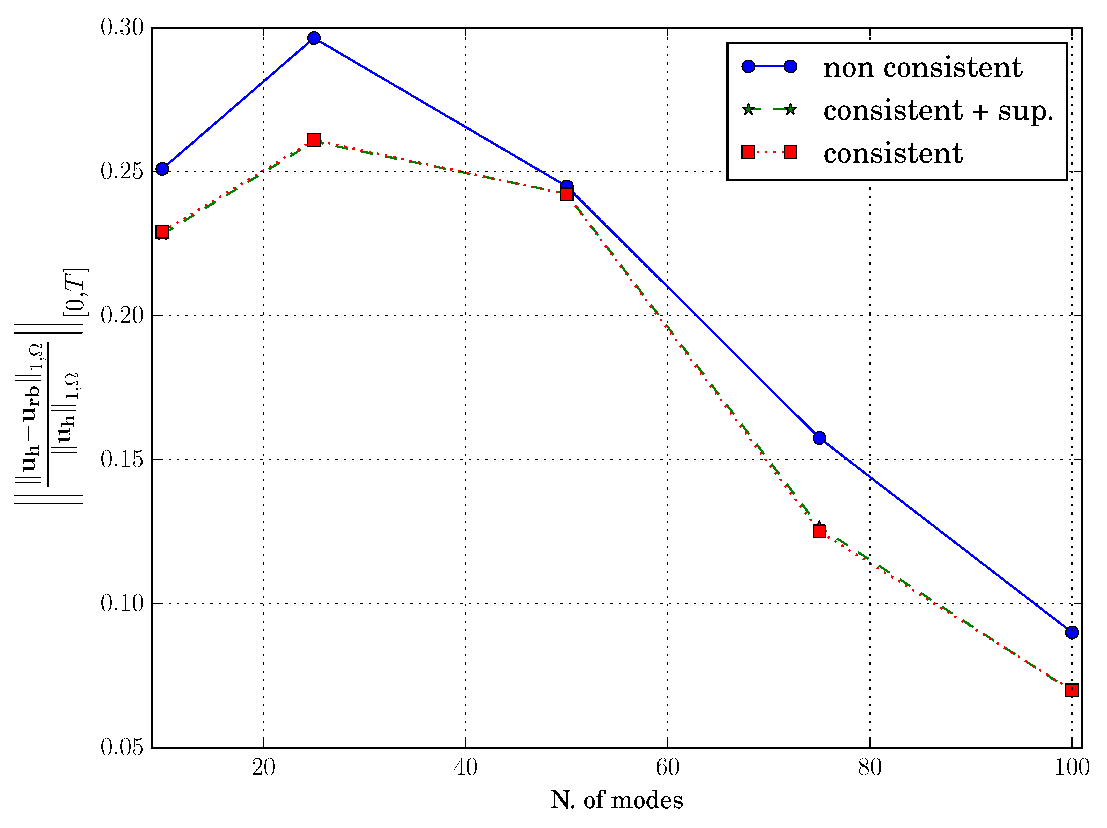}
\end{minipage}
\begin{minipage}{0.5\textwidth}
\centering
\includegraphics[width=\textwidth]{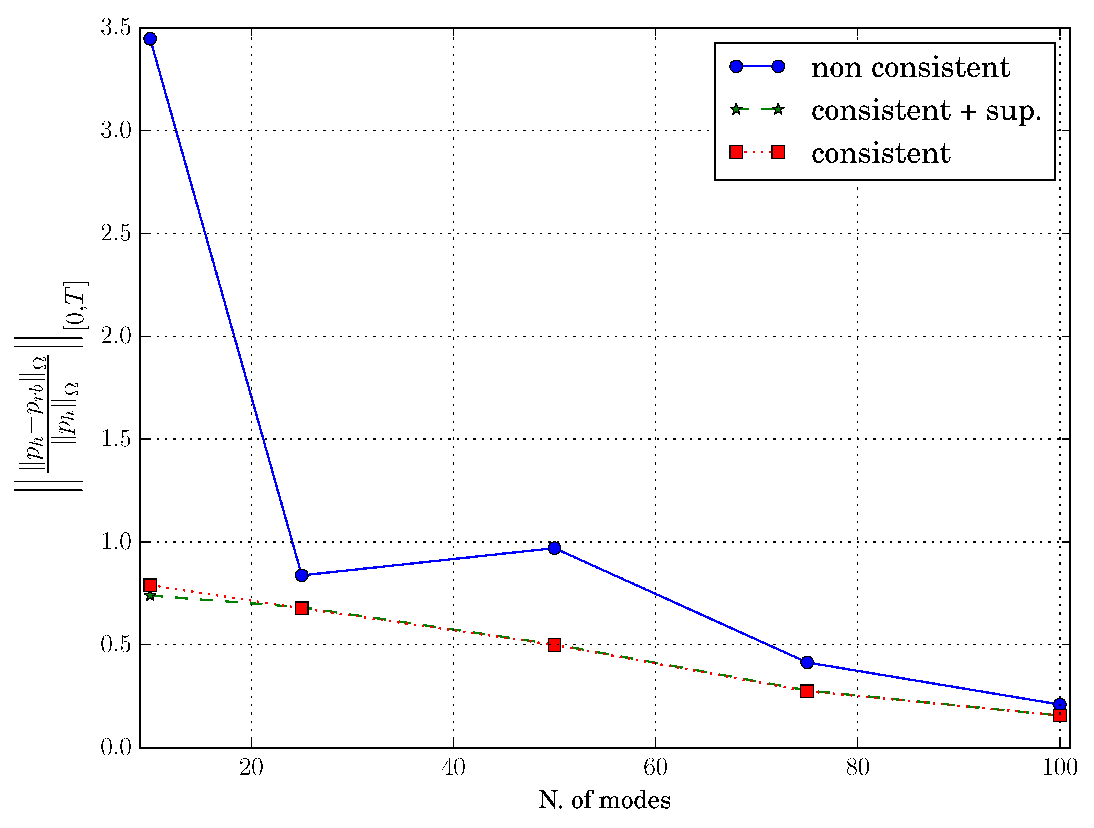}
\end{minipage}
\end{minipage}
\caption{Velocity relative error in the $L^2([0, T], \gr{H}^1(\Omega))$ norm (left) and pressure relative error in the $L^2([0, T], L^2(\Omega))$ norm (right), as a function of the number $N$ of POD modes. }
\label{fig:error_modes}
\end{figure}

\begin{figure}
\begin{minipage}{\textwidth}
\begin{minipage}{0.5\textwidth}
\centering
\includegraphics[width=\textwidth]{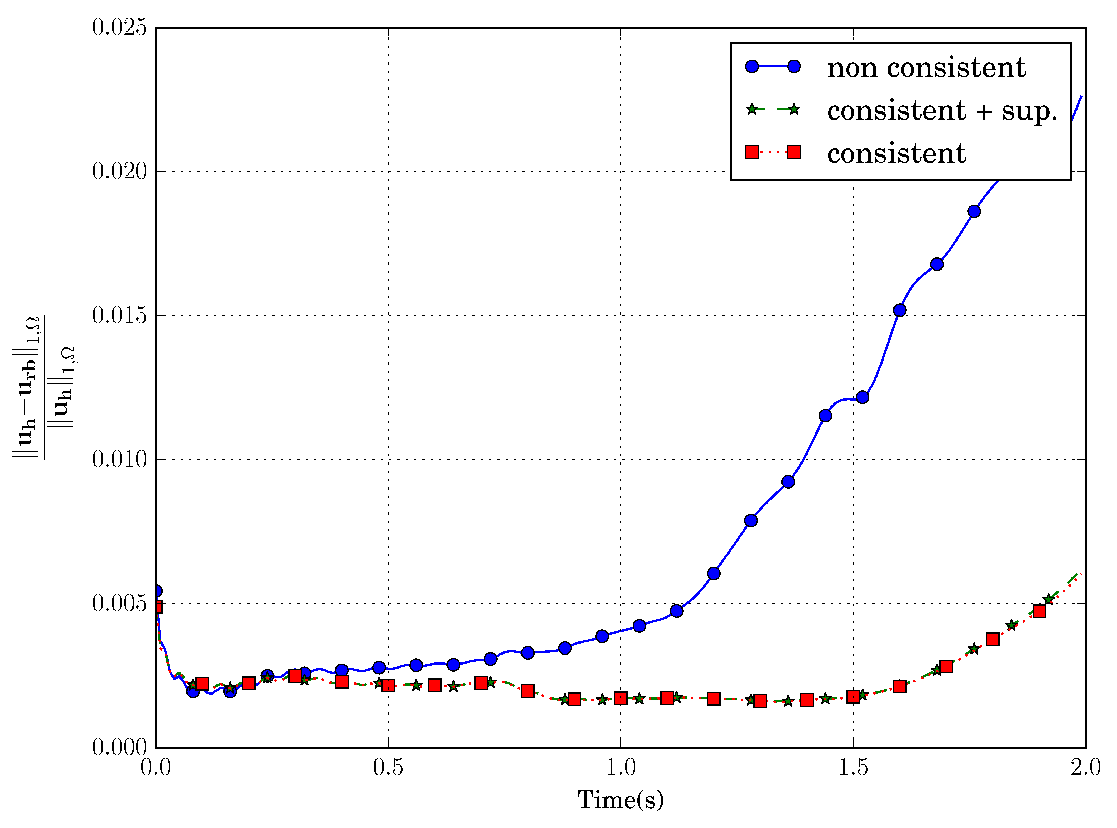}
\end{minipage}
\begin{minipage}{0.5\textwidth}
\centering
\includegraphics[width=\textwidth]{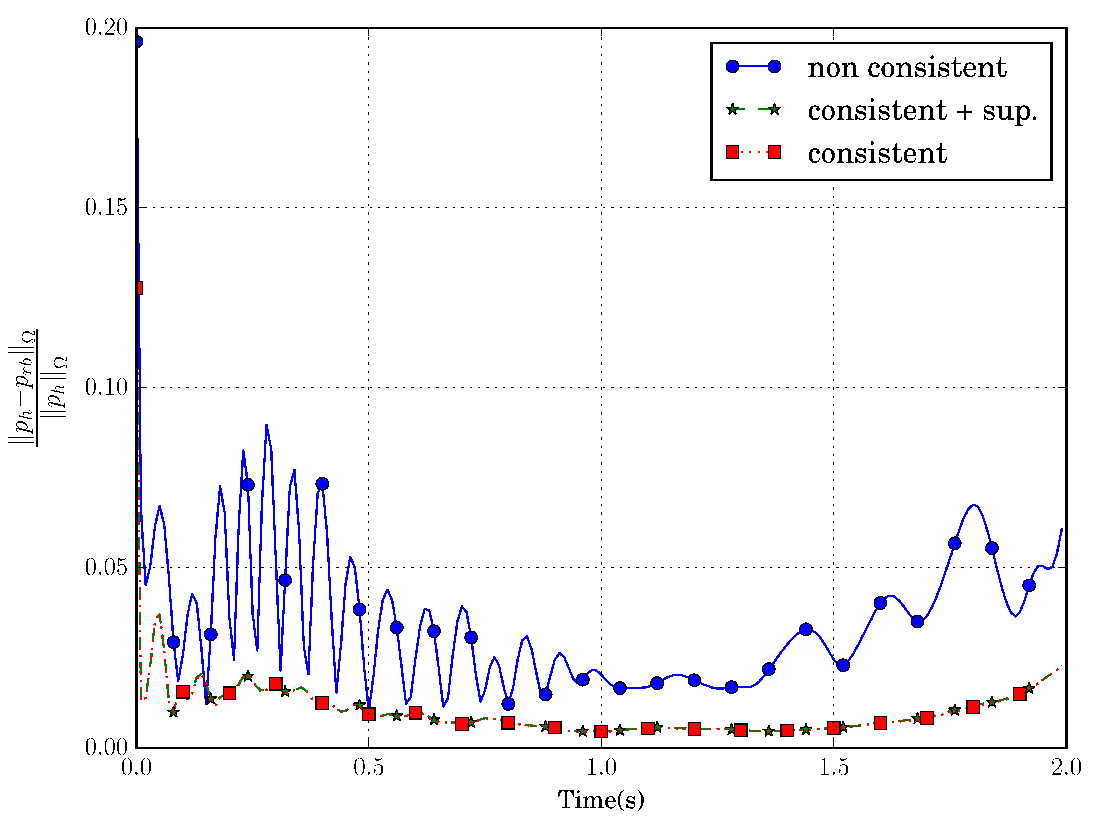}
\end{minipage}
\end{minipage}
\caption{Velocity relative error in the $\gr{H}^1(\Omega)$ norm (left) and pressure relative error in the $L^2(\Omega)$ norm (right), as a function of time. ROM evaluation was carried out for $N=100$. }
\label{fig:error_time}
\end{figure}

\begin{figure}
\begin{minipage}{\textwidth}
\begin{minipage}{0.5\textwidth}
\centering
\includegraphics[width=\textwidth]{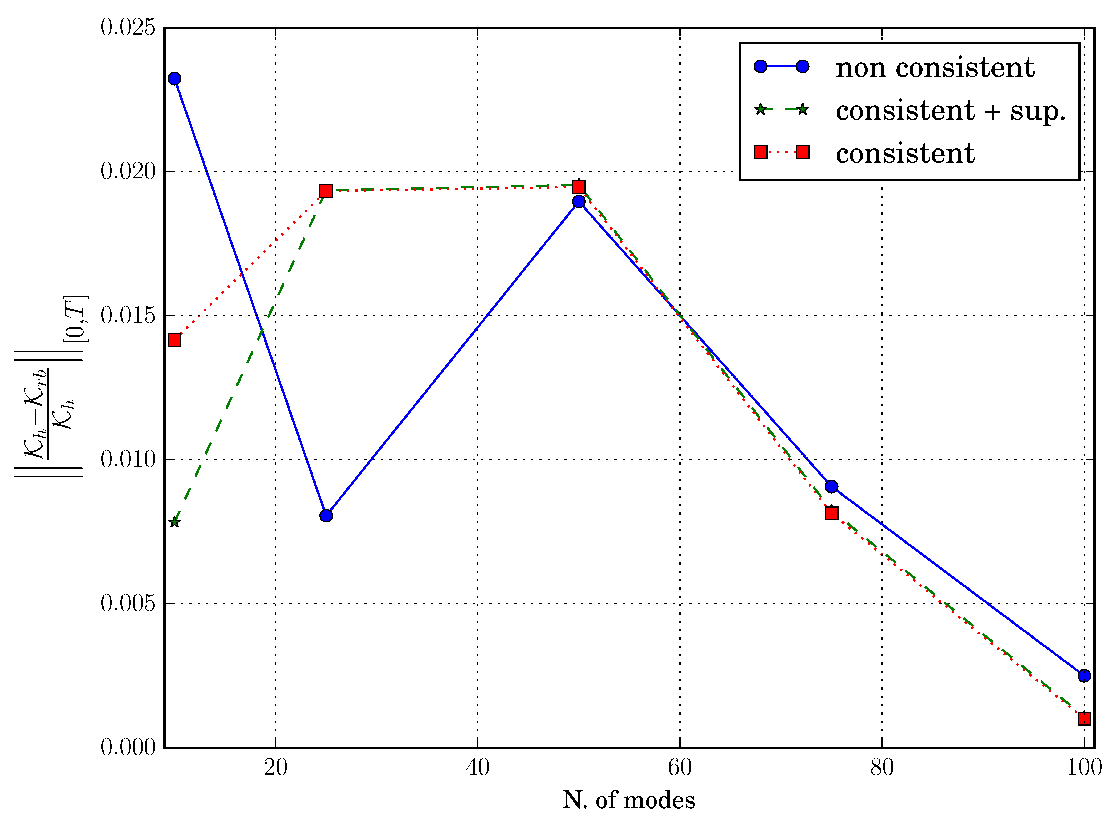}
\end{minipage}
\begin{minipage}{0.5\textwidth}
\centering
\includegraphics[width=\textwidth]{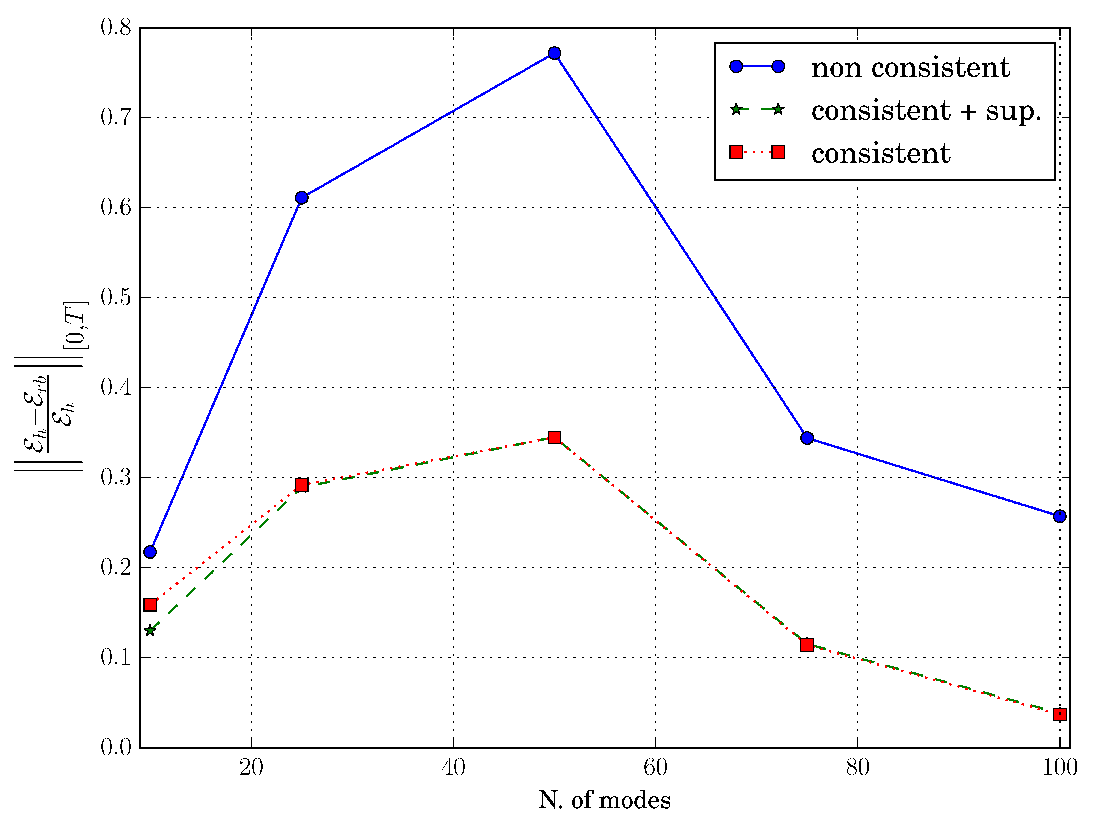}
\end{minipage}
\end{minipage}
\caption{Kinetic energy (left) and enstrophy (right) relative errors as a function of the number of employed modes.}
\label{fig:error_kinetic_enstophy_modes}
\end{figure}

\begin{figure}
\begin{minipage}{\textwidth}
\begin{minipage}{0.5\textwidth}
\centering
\includegraphics[width=\textwidth]{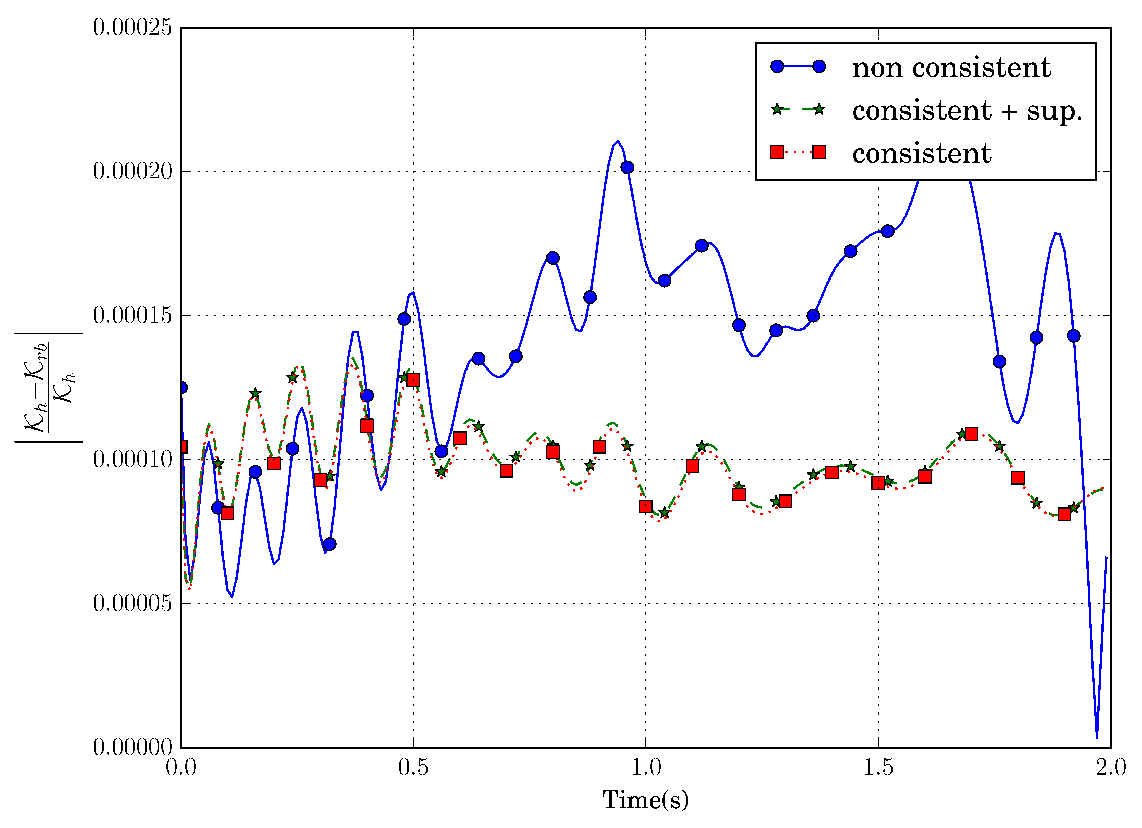}
\end{minipage}
\begin{minipage}{0.5\textwidth}
\centering
\includegraphics[width=\textwidth]{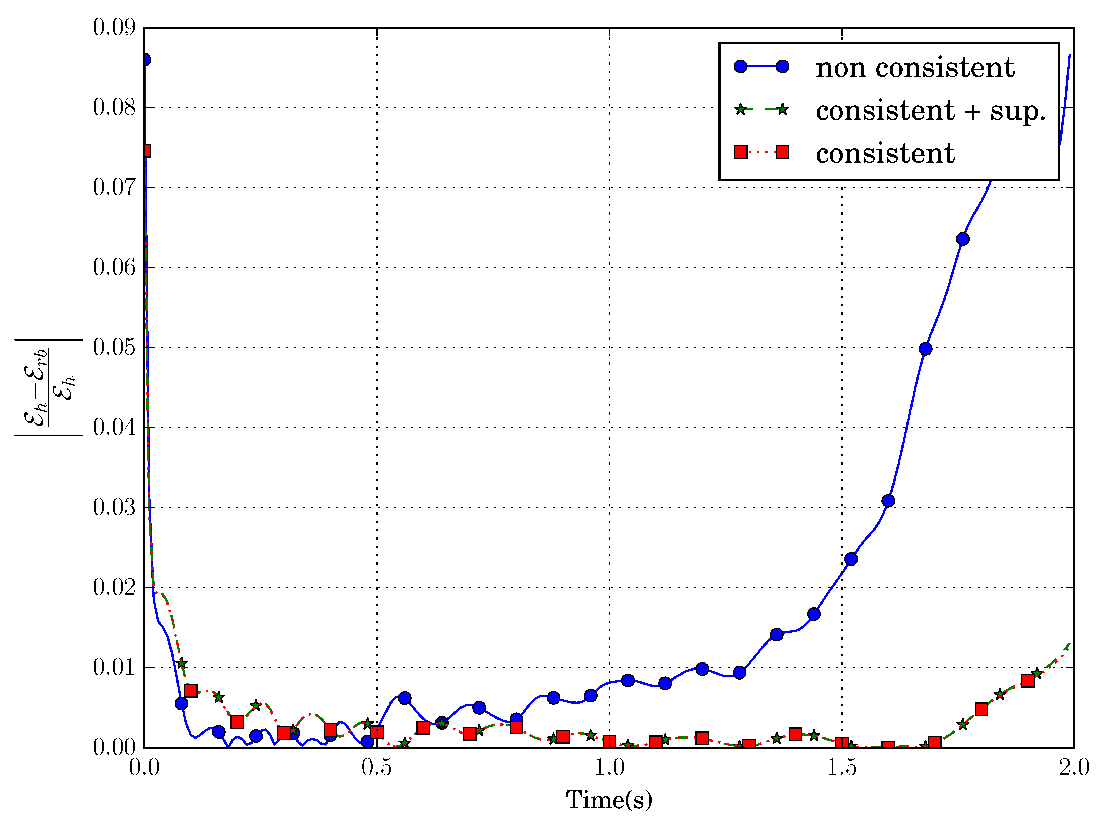}
\end{minipage}
\end{minipage}
\caption{Kinetic energy (left) and enstrophy (right) relative errors as a function of time. ROM evaluation was carried out for $N=100$.}
\label{fig:error_kinetic_enstophy}
\end{figure}

The evaluation of the ROM is then queried. Relative errors for various quantities of interest are depicted Figures \ref{fig:error_modes}-\ref{fig:error_kinetic_enstophy}. Three possible ROMs are compared, corresponding to (i) consistent ROM with supremizers, (ii) consistent ROM without supremizers, (iii) non-consistent ROM with supremizers. The fourth option, i.e. non-consistent ROM without supremizer, is omitted because of its lack of pressure stabilization, which would result in a singular ROM. More in detail, Figure \ref{fig:error_modes} shows the velocity and pressure relative errors, integrated over time, for increasing number of modes $N$ from 10 to 100. 

Furthermore, for fixed $N=100$, Figure \ref{fig:error_time} plots the velocity and pressure relative errors for each time $t \in [0, 2]$, while Figure \ref{fig:error_kinetic_enstophy} carries out a similar analysis for the kinetic energy $\mathcal{K}(\gr{u}) = \int_\Omega \lvert \gr{u} \rvert^2$ and enstrophy $\mathcal{E}(\gr{u}) = \int_\Omega \lvert \text{curl}(\gr{u})\rvert^2$. 

Results in Figures \ref{fig:error_modes}-\ref{fig:error_kinetic_enstophy} show that the consistent ROM is {systematically able to reproduce the FOM results with an higher accuracy} than the non-consistent one. In particular, improvements of {$73.65$\%, $65.12$\%, $2.83$\%, $85.48$\%} can be observed at the final time for the consistent ROM {(with $N=100$)} over the non-consistent one for velocity, pressure, kinetic energy and enstrophy, respectively. Thus, as in \cite{Shafqat_unsteady,Shafqat,PACCIARINI20141}, it is our conclusion that ensuring consistency between the full order and reduced order formulations results {in a FOM characterized by the best accuracy with respect to the FOM}.

Figures \ref{fig:error_modes}-\ref{fig:error_kinetic_enstophy} also highlight that the addition of supremizers is not necessary once a consistent ROM is sought. Indeed, lines corresponding to the options (i) consistent ROM with supremizer, and (ii) consistent ROM without supremizers, are always overlapping. A careful analysis of the coefficients of the reduced velocity in case (i) reveals that degrees of freedom corresponding to supremizer modes are numerically zero, thus not contributing to the solution. We claim that this is thanks to the stabilization through strong residual operators in the VMS formulation, which (through Galerkin projection) results in an inf-sup stable ROM formulation.\footnote{Nonetheless, we remark that this topic deserves further investigation. Indeed, in case of geometrical parametrization with a similar residual based SUPG stabilization we observe non-zero supremizer coefficients \cite{Shafqat}. We do mention that, in our experience, ROMs seem to be particularly sensitive to the choice of the stabilization at the full order level, especially for what concerns their consistency (or lack thereof) with the strong form. For instance, in \cite{KaBaRo18} we pursue a weaker (non-residual based, non-fully consistent) stabilization in the full order discretization, and we observe that the Galerkin projection of such stabilization terms does not result in a stable ROM, yielding a mandatory supremizer enrichment. In \cite{KaStaNouScoRo2018}, where a stabilized Shifted Boundary FEM is used as full order discretization in a Stokes setting (see also \cite{KaStaNoRoSco18_Heat} for a Poisson setting), we observe that the supremizer enrichment is not strictly needed but it leads to a better approximation of the pressure field.}\label{foot:sup} {A more detailed numerical analysis on this topic is subject of on-going work and will be addressed in a future paper.}

Finally, Figures \ref{fig:plots_u} and \ref{fig:plots_p} show a comparison of the flow fields for $t=1.0$ s, $t=1.5$ s, $t=2.0$ for full order model, consistent ROM with supremizer, non-consistent ROM with supremizer; the remaining viable option (consistent ROM without supremizer) is omitted as the difference among consistent ROMs is negligible. The plots show that, especially for large times (see e.g. $t=2$ s), flow fields produced by the non-consistent ROM are qualitatively different from the full order simulation, thus resulting in larger relative errors.

\section{Conclusions and future developments}\label{sec:concl}
In the present work we introduced a novel reduced order model dealing with moderately high Reynolds numbers. The ROM has been constructed starting from VMS stabilized full order simulations and different strategies have been tested during the projection stage. 

The consistent strategy, which is based on the projection of also the full order stabilization terms onto the reduced basis spaces, proved to be the best choice in terms of accuracy of the results {for what concerns the capability of the ROM to reproduce the solution obtained by the FOM model}.

The consistent strategy permits also to avoid the supremizer enrichment of the reduced velocity space. We believe that this fact is justified by the VMS full order discretization which, being inf-sup stabilized, transfers this property also to the reduced order model. However, this topic deserves still further investigation especially in the case of geometrical parametrization (see footnote 6 for further details).  

As previously mentioned, as future perspective our interest is to extend the current approach also to geometrically parametrized problems and to investigate the efficiency and applicability of hyper reduction techniques. {Further future perspectives are related to the investigation of different closure models (as well as comparison with DNS results, when possible) and 3D simulations.}

\section*{Acknowledgements}
We acknowledge Prof. Guglielmo Scovazzi from Duke University for the fruitful discussions concerning the implementation of the VMS method. We acknowledge the support by European Union Funding for Research and Innovation - Horizon 2020 Program - in the framework of European Research Council Executive Agency: H2020 ERC Consolidator Grant 2015 AROMA-CFD project 681447 ``Advanced Reduced Order Methods with Applications in Computational Fluid Dynamics''. We also acknowledge the INDAM-GNCS projects ``Metodi numerici avanzati combinati con tecniche di riduzione computazionale per PDEs parametrizzate e applicazioni'' and ``Numerical methods for model order reduction of PDEs''. The computations in this work have been performed with RBniCS \cite{rbnics} library, developed at SISSA mathLab, which is an implementation in FEniCS \cite{logg2012automated} of several reduced order modelling techniques; we acknowledge developers and contributors to both libraries. 

\clearpage
\bibliographystyle{abbrv}
\bibliography{bibfile}

\end{document}